\documentclass[11pt]{amsart}
\usepackage{graphicx,amssymb,amsthm}
\usepackage{amscd}
\usepackage[matrix,arrow,curve]{xy}
\usepackage{a4wide}

\def\baselinestretch{1.15}

\theoremstyle{plain}
\newtheorem{theorem}{Theorem}[section]
\newtheorem{lemma}[theorem]{Lemma}
\newtheorem{proposition}[theorem]{Proposition}
\newtheorem{corollary}[theorem]{Corollary}
\theoremstyle{definition}
\newtheorem{definition}[theorem]{Definition}
\newtheorem{example}[theorem]{Example}
\newenvironment{mytheorem}[1]
    {\par\addvspace{\smallskipamount}\noindent\textbf{#1.}\ \em\ignorespaces}
    {\par\addvspace{\smallskipamount}}

\def\C{\mathcal C}
\def\D{\mathcal D}
\def\Z{\mathbb Z}

\def\c{\mathbf c}
\def\d{\mathbf d}
\def\n{\mathbf n}

\def\infs{\inf{\!}_s}
\def\sups{\sup{\!}_s}
\def\myangle#1{\langle #1\rangle}
\def\ext{{\operatorname{ext}}}
\def\R{\operatorname{\mathcal R}}
\def\Rext{\R_\ext}
\def\LM{\operatorname{s}_{\!L}}
\def\RM{\operatorname{s}_{\!R}}

\def\wedgeL{\wedge_L}
\def\wedgeR{\wedge_R}
\def\veeL{\vee_{\!L}}
\def\veeR{\vee_{\!R}}
\let\le\leqslant
\let\ge\geqslant

\def\St{{\operatorname{St}}}
\def\ST{{\operatorname{\it St}}}
\def\ST{{St}}
\def\brindex{\operatorname{br}}

\def\t{\operatorname{\it t}}
\def\INF{\t_{\inf}}
\def\SUP{\t_{\sup}}

\begin{document}

\title
[A Garside-theoretic approach to the reducibility problem in braid groups]
{A Garside-theoretic approach\\ to the reducibility problem in braid groups}

\author{Eon-Kyung Lee and Sang-Jin Lee}

\address{Department of Applied Mathematics, Sejong University,
Seoul 143-747, Korea}
\email{eonkyung@sejong.ac.kr}

\address{Department of Mathematics, Konkuk University,
Seoul 143-701, Korea}
\email{sangjin@konkuk.ac.kr}

\date{\today}

\begin{abstract}
Let $D_n$ denote the $n$-punctured disk in the complex plane,
where the punctures are on the real axis.
An $n$-braid $\alpha$ is said to be \emph{reducible} if
there exists an essential curve system $\C$ in $D_n$,
called a \emph{reduction system} of $\alpha$,
such that $\alpha*\C=\C$
where $\alpha*\C$ denotes the action of the braid $\alpha$
on the curve system $\C$.
A curve system $\C$ in $D_n$ is said to be \emph{standard}
if each of its components is isotopic to a round circle
centered at the real axis.

In this paper, we study the characteristics of
the braids sending a curve system to a standard curve system,
and then the characteristics of the conjugacy classes
of reducible braids.
For an essential curve system $\C$ in $D_n$,
we define the \emph{standardizer} of $\C$ as
$\St(\C)=\{P\in B_n^+:P*\C\mbox{ is standard}\}$ and
show that $\St(\C)$ is a sublattice of $B_n^+$.
In particular, there exists a unique minimal element in $\St(\C)$.
Exploiting the minimal elements of standardizers
together with canonical reduction systems
of reducible braids,
we define the outermost component of reducible braids,
and then show that, for the reducible braids
whose outermost component is simpler than the whole braid
(including split braids),
each element of its ultra summit set has a standard reduction system.
This implies that, for such braids,
finding a reduction system is
as easy as finding a single element of the ultra summit set.

\end{abstract}

\maketitle

\section{Introduction}
Let $D_n=\{z\in\mathbb C: |z|\le n+1\}\setminus\{1,\ldots,n\}$,
the $n$-punctured disk in the complex plane with punctures
lying on the real axis.
The $n$-braid group $B_n$ acts on the set of curve systems in $D_n$.
For an $n$-braid $\alpha$ and a curve system $\C$ in $D_n$,
let $\alpha*\C$ denote the action of $\alpha$ on $\C$.
An $n$-braid $\alpha$ is said to be \emph{reducible}
if $\alpha*\C=\C$ for some essential curve system $\C$ in $D_n$,
called a \emph{reduction system} of $\alpha$.
In this paper, we are interested in the \emph{reducibility problem}:
given a braid, decide whether it is reducible or not and
find a reduction system if it is reducible.

\subsection{Motivation and some of previous works}
The Nielsen-Thurston classification
theorem~\cite{Thu88} states that an irreducible automorphism of
an orientable surface with negative euler characteristic
is either periodic or pseudo-Anosov up to isotopy. Recall that
an orientation preserving self-diffeomorphism $f$ of a surface $S$ is said to be
\begin{itemize}
\item
\emph{periodic} if $f^k$ is isotopic to the identity for some $k\ne 0$;

\item
\emph{reducible} if there exist pairwise disjoint simple closed curves
$C_1,\ldots,C_k$ in $S$, isotopic to neither a point nor a puncture
nor a boundary component,
such that $f(\C)$ is isotopic to $\C$,
where $\C=C_1\cup\cdots \cup C_k$;

\item
\emph{pseudo-Anosov} if there exist a pair of transverse measured foliations
$(F^s,\mu^s)$ and $(F^u,\mu^u)$ and a real $\lambda>1$ such that
$f(F^s,\mu^s)=(F^s,\lambda^{-1}\mu^s)$ and
$f(F^u,\mu^u)=(F^u,\lambda\mu^u)$.
\end{itemize}

There have been several approaches to the problem of deciding
dynamical types of surface automorphisms.
Bestvina and Handel~\cite{BH95} made the train track algorithm that,
given a surface automorphism,
decides its dynamical type and finds its dynamical
structure: a pair of transverse measured foliations for
a pseudo-Anosov automorphism; a reduction system for a reducible
automorphism.
Benardete, Guti\'errez and Nitecki~\cite{BGN95}
solved the reducibility problem in braid groups.
(It is known that a periodic $n$-braid is
conjugate to either $(\sigma_1\sigma_2\cdots \sigma_{n-1})^l$
or $(\sigma_1(\sigma_1\sigma_2\cdots\sigma_{n-1}))^l$ for some
integer $l$~\cite{Ker19,Eil34,BDM02}.
This implies that $\alpha$ is a periodic $n$-braid
if and only if either $\alpha^n$ or $\alpha^{n-1}$
is equal to $\Delta^{2m}$ for some integer $m$.
Hence, it is easy to decide the periodicity of braids.
Therefore, in order to decide the dynamical type of a given braid,
it suffices to decide the reducibility.)
Humphries~\cite{Hum91} solved the problem of recognizing
split braids.

With the above results,
solving the reducibility problem and the problem of recognizing split braids
seems at least as hard as solving the conjugacy problem.
When using the train track algorithm,
one needs to describe a given $n$-braid as a graph map of the $n$-bouquet,
and the length of this description grows exponentially with respect
to the length of the braid word on Artin generators.
The other two solutions need to use the algorithms solving the conjugacy
problem in braid groups.

\smallskip
Another motivation for this work is the close relationship between
the reducibility problem and the conjugacy problem.
The approach to the conjugacy problem in braid groups
can be divided into two steps:
solving the reducibility problem and solving the conjugacy problem
for irreducible braids.
See~\cite[\S1.4]{BGG06a}
for a more precise description of this strategy.
The conjugacy problem for periodic braids is easy to solve.
There are two different polynomial-time solutions to this case
by Birman, Gebhardt and Gonz\'alez-Meneses~\cite{BGG06c}
and by the authors~\cite{LL07b}.
For the conjugacy problem for pseudo-Anosov mapping classes,
there are several results.
In~\cite{Los93}, Los solved the problem for pseudo-Anosov braids
by using combinatorial efficient representatives.
Recently, Fehrenbach and Los~\cite{FL07} proposed an algorithm that finds
roots and symmetries of pseudo-Anosov mapping classes
together with a new solution to the conjugacy problem.
Mazur and Minsky~\cite{MM99,MM00} showed that, fixing a
mapping class group and a finite set of generators for this group,
there exists a constant $K$ such that
if $\alpha$ and $\beta$ are conjugate pseudo-Anosov mapping classes
then there is a conjugating element $\gamma$ with $|\gamma|\le K
(|\alpha|+|\beta|)$, where $|\cdot|$ denotes the word length.
In order to extend the results on irreducible braids to general braids,
we need to solve the reducibility problem more efficiently.

\smallskip
For the last ten years, no serious progress has been made
in the reducibility problem.
On the other hand, recently, there have been several new contributions to Garside-theoretic
approach to braid groups, for example~\cite{Deh02, FG03, Geb05, Lee07}.
Exploiting them, we study the characteristics of the conjugacy classes
of reducible braids.
Our approach uses neither the train track algorithm
nor the complete conjugacy algorithm.
We hope that our results are useful in obtaining
a more efficient solution to the reducibility problem in braid groups.

\subsection{Our results}
Before stating our results,
we recall some notions and results from the Garside theory in braid groups.

\begin{itemize}
\item
Let $B_n^+$ be the submonoid of $B_n$ generated by
$\sigma_1,\ldots,\sigma_{n-1}$.
The partial order $\le_R$ on $B_n^+$ is defined as follows:
for $P,Q\in B_n^+$, $P\le_R Q$ if $Q=SP$ for some $S\in B_n^+$.
The poset $(B_n^+,\le_R)$ is a lattice,
i.e., there exist the gcd $P\wedgeR Q$ and the lcm
$P\veeR Q$ of $P,Q\in B_n^+$.

\item
For $\alpha\in B_n$, there are integer-valued invariants $\inf(\alpha)$
and $\sup(\alpha)$.
Let $[\alpha]$ denote
the conjugacy class of $\alpha\in B_n$.
The following are conjugacy invariants.
$$
\begin{array}{ll}
\infs(\alpha)=\max\{\inf(\beta):\beta\in[\alpha]\}\qquad
&\INF(\alpha)=\lim_{m\to\infty} \inf(\alpha^m)/m\\
\sups(\alpha)=\min\{\sup(\beta):\beta\in[\alpha]\}
&\SUP(\alpha)=\lim_{m\to\infty} \sup(\alpha^m)/m
\end{array}
$$

\item
In the conjugacy class $[\alpha]$, there are finite, nonempty,
computable subsets,
the super summit set $[\alpha]^S$,
the ultra summit set $[\alpha]^U$
and the stable super summit set $[\alpha]^\ST$.
They depend only on the conjugacy class,
and $[\alpha]^U, [\alpha]^\ST\subset[\alpha]^S$.

\end{itemize}

\begin{figure}
\includegraphics[scale=.6]{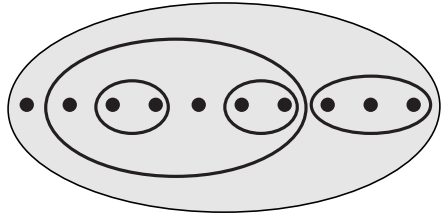}
\vskip -2mm
\caption{A standard curve system in $D_{10}$}
\label{fig:standard}
\end{figure}

We call an essential curve system (see Definition~\ref{def:curves}) in $D_n$
a \emph{standard} curve system
if each component is isotopic to a round circle centered
at the real axis as in Figure~\ref{fig:standard}.
For an essential curve system $\C$ in $D_n$,
we define the \emph{standardizer} of $\C$ as the set
$$
\St(\C)= \{P\in B_n^+:\mbox{$P*\C$ is standard}\}
$$
where $P*\C$ denotes the left action of the positive braid $P$
on the curve system $\C$, and then show the following.

\begin{mytheorem}{Theorem~\ref{thm:closed}}
For an essential curve system $\C$ in $D_n$,
its standardizer $\St(\C)$ is closed under $\wedgeR$ and $\veeR$,
and hence a sublattice of $B_n^+$.
Therefore $\St(\C)$ contains a unique $\le_R$-minimal element.
\end{mytheorem}

\begin{mytheorem}{Theorem~\ref{thm:StMain}}
Let $\alpha$ be a reducible $n$-braid with a reduction system $\C$.
Let $P$ be the $\le_R$-minimal element of\/ $\St(\C)$.
Then the following hold.
\begin{enumerate}
\item[(i)] $\inf(\alpha)\le\inf(P\alpha P^{-1})\le\sup(P\alpha P^{-1})\le\sup(\alpha)$.
\item[(ii)] If\/ $\alpha\in[\alpha]^S$, then $P\alpha P^{-1}\in[\alpha]^S$.
\item[(iii)] If\/ $\alpha\in[\alpha]^U$, then $P\alpha P^{-1}\in[\alpha]^U$.
\item[(iv)] If\/ $\alpha\in[\alpha]^\ST$, then $P\alpha P^{-1}\in[\alpha]^\ST$.
\end{enumerate}
\end{mytheorem}

Theorem~\ref{thm:closed} is essential in our approach to the reducibility problem,
as the closedness under $\wedgeR$ of
$\{P\in B_n^+: P\beta P^{-1}\in[\alpha]^S\}$ and
$\{P\in B_n^+: P\beta P^{-1}\in[\alpha]^U\}$ for $\beta\in[\alpha]^S$
plays an important role in solving the conjugacy problem~\cite{FG03, Geb05}.
Theorem~\ref{thm:StMain} shows that
standardizing a reduction system $\C$ of a braid
by the $\le_R$-minimal element of $\St(\C)$
preserves the membership of the super summit set,
ultra summit set and stable super summit set.

\smallskip
It is known by Birman, Lubotzky and McCarthy~\cite{BLM83}
and Ivanov~\cite{Iva92} that a reducible surface automorphism
admits a unique \emph{canonical reduction system}.
For $\alpha\in B_n$, let $\Rext(\alpha)$ be the collection of the outermost
components of the canonical reduction system of $\alpha$.
Let $P$ be the $\le_R$-minimal element of $\St(\Rext(\alpha))$.
Since $\Rext(P\alpha P^{-1})=P*\Rext(\alpha)$ is standard,
the outermost component of $D_n\setminus\Rext(P\alpha P^{-1})$
is naturally identified with the $k$-punctured disk $D_k$ for some $k\le n$.
We define the \emph{outermost component} $\alpha_\ext$ of $\alpha$
as the $k$-braid obtained by restricting the braid $P\alpha P^{-1}$ to
the outermost component of $D_n\setminus \Rext(P\alpha P^{-1})$.
See \S5 for the precise definition.
The following is the main result of this paper.
(In the statement, $[\alpha]^U_\d$ denotes
the ultra summit set of $\alpha$ with respect to decycling.
See the next section for the precise definition.)

\begin{mytheorem}{Theorem~\ref{thm:main}}
Let $\alpha$ be a non-periodic reducible $n$-braid.
\begin{enumerate}
\item[(i)]
If\/ $\infs(\alpha_\ext)>\infs(\alpha)$,
then each element of\/ $[\alpha]^U$ has a standard reduction system.

\item[(ii)]
If\/ $\sups(\alpha_\ext)<\sups(\alpha)$,
then each element of\/ $[\alpha]_\d^U$ has a standard reduction system.

\item[(iii)]
If\/ $\alpha$ is a split braid,
then each element of\/ $[\alpha]^U \cup [\alpha]_\d^U$
has a standard reduction system.

\item[(iv)]
If\/ $\alpha_\ext$ is periodic,
then there exists $1\le q< n$ such that each element
of\/ $[\alpha^q]^U \cup [\alpha^q]^U_\d$ has a standard reduction system.

\item[(v)]
If\/ $\INF(\alpha_\ext)>\INF(\alpha)$,
then there exists $1\le q< n(n-1)/2$
such that each element
of\/ $[\alpha^q]^U$ has a standard reduction system.

\item[(vi)]
If\/ $\SUP(\alpha_\ext)<\SUP(\alpha)$,
then there exists $1\le q< n(n-1)/2$
such that each element of\/ $[\alpha^q]^U_\d$
has a standard reduction system.
\end{enumerate}
\end{mytheorem}

Roughly speaking, the first statement of the above theorem says that
if the outermost component $\alpha_\ext$ is simpler than
the whole braid $\alpha$ up to conjugacy
from a Garside-theoretic point of view,
then every element of $[\alpha]^U$ has a standard reduction system.
In this case, finding a reduction system is
as easy as finding one element in the ultra summit set,
because it is easy to find a standard reduction system
of a given braid if it exists by the results in~\cite{BGN93}.
In \S7, we present three examples showing that
the conditions in Theorem~\ref{thm:main} cannot be weakened.

In~\cite{BGN95}, Benardete, Guti\'errez and Nitecki
showed that \emph{if a braid is reducible,
then there exists an element in its super summit set
which has a standard reduction system.}
(The notion of ultra summit set appeared later than their work,
and from their proof we can replace `super summit set'
in their statement with `ultra summit set'.)
While their result concerns the \emph{existence} of an ultra summit
element with a standard reduction system,
Theorem~\ref{thm:main}~(i)-(iii)
show that, under a certain condition,
\emph{every} ultra summit element
has a standard reduction system.

We remark that the six types of braids in Theorem~\ref{thm:main} cover
most reducible braids.
The braid $\alpha_\ext$ can be obtained, up to conjugacy,
by deleting some strands from $\alpha$,
hence $\alpha_\ext$ cannot be more complicated than $\alpha$.
Indeed, the following inequalities always hold
(see Lemma~\ref{thm:infscompare}):
$$
\begin{array}{cc}
\infs(\alpha_\ext)\ge\infs(\alpha);
&\sups(\alpha_\ext)\le\sups(\alpha);\\
\INF(\alpha_\ext)\ge \INF(\alpha);
& \SUP(\alpha_\ext)\le\SUP(\alpha).
\end{array}
$$
Theorem~\ref{thm:main} shows the characteristics of the braid conjugacy classes
for which at least one
of the above inequalities is strict.

\smallskip
We briefly explain the idea of proof of Theorem~\ref{thm:main}.

\begin{itemize}

\item In \S6, we show that if $\alpha$ is a split braid with
the minimal word length in the conjugacy class,
then the outermost component $\R_\ext(\alpha)$ of the canonical
reduction system of $\alpha$ is standard.
Since a positive braid has the minimal word length in the conjugacy class,
we have the following: \emph{if $P$ is a positive split braid,
then $\R_\ext(P)$ is standard.}

\item If a braid $\alpha$ commutes with a non-periodic reducible braid $\beta$,
then the canonical reduction system of $\beta$ is a reduction system of $\alpha$.
Combining this with the previous observation, we have the following:
\emph{if $\alpha P=P\alpha$ for some positive split braid $P$,
then $\Rext(P)$ is a standard reduction system of $\alpha$.}

\item If $\alpha$ belongs to the ultra summit set,
then there exists a finite sequence
$\alpha=\alpha_0\to\alpha_1\to\cdots\to\alpha_m=\alpha$
for some $m\ge 1$, where $\alpha_{i+1}=A_i\alpha_iA_i^{-1}$
for some permutation braid $A_i$
for $i=0,\ldots,m-1$.
If we let $T=A_{m-1}\cdots A_1A_0$, then $T\alpha=\alpha T$.
Exploiting the $\le_R$-minimal elements of the standardizers
$\St(\Rext(\alpha_i))$,
we show that $T$ is a positive split braid
if $\infs(\alpha_\ext)>\infs(\alpha)$,
from which Theorem~\ref{thm:main}~(i) follows.
The other statements are proved using this.
\end{itemize}

\subsection{Organization}
In \S2, we review the Garside theory in brad groups.
In \S3, we study the normal form of the braids that send a standard curve system
to a standard curve system.
In \S4, we prove Theorems~\ref{thm:closed} and~\ref{thm:StMain}.
In \S5, we study the properties of the outermost component $\alpha_\ext$
of a non-periodic reducible braid $\alpha$.
In \S6, we show that if a split braid has the minimal
word length in the conjugacy class, then the outermost component of its
canonical reduction system is standard.
In \S7 and \S8, we prove Theorem~\ref{thm:main}, using the results
of the previous sections.

\subsection*{Acknowledgements}

We are most grateful to the anonymous referee of this journal
for valuable comments and suggestions on the paper,
especially for pointing out that
our initial proof of Theorem~\ref{thm:StMain} contains a mistake.
The proof is corrected as suggested by the referee.
We are also very thankful to Won Taek Song for helpful
conversations, and Juan Gonz\'alez-Meneses and Bert Wiest for
providing Example~\ref{ex:3}.
This work was supported by the Korea Science and Engineering
Foundation (KOSEF) grant funded by the Korea government (MOST)
(No.~R01-2007-000-20293-0).

\section{Garside theory in braid groups}
We give necessary definitions and results on Garside theory
in braid groups.
See~\cite{Gar69, Thu92, EM94, BKL98, DP99, Deh02, FG03, Geb05} for details.
The $n$-braid group $B_n$ has the group presentation
$$
B_n  =  \left\langle \sigma_1,\ldots,\sigma_{n-1} \left|
\begin{array}{ll}
\sigma_i \sigma_j = \sigma_j \sigma_i & \mbox{if } |i-j| \ge 2, \\
\sigma_i \sigma_j \sigma_i = \sigma_j \sigma_i \sigma_j & \mbox{if } |i-j| = 1.
\end{array}
\right.\right\rangle,
$$
where $\sigma_i$ is the isotopy class of the positive half Dehn-twist
along the straight line segment connecting the punctures $i$ and
$i+1$. An $n$-braid can be regarded as a collection of $n$ strands
$l=l_1\cup\cdots\cup l_n$ in $[0,1]\times D^2$
such that $|\,l\cap (\{t\}\times D^2)|=n$ for $0\le t\le 1$ and
$l\cap (\{0,1\}\times D^2)=\{0,1\}\times \{1,\ldots,n\}$.

\subsection{Positive braid monoid}
Let $B_n^+$ be the monoid generated by $\sigma_1,\ldots,\sigma_{n-1}$
with the defining relations: $\sigma_i\sigma_j=\sigma_j\sigma_i$ for $|i-j|\ge 2$;
$\sigma_i\sigma_j\sigma_i=\sigma_j\sigma_i\sigma_j$ for $|i-j|=1$.
$B_n^+$ is a (left and right) cancellative monoid
that embeds in $B_n$ under the canonical homomorphism.
$B_n^+$ is called the \emph{positive braid monoid}
and its elements are called \emph{positive braids}.

\begin{definition}
The partial orders $\le_L$ and $\le_R$ on $B_n^+$ are defined as follows:
for $P,Q\in B_n^+$,
$P\le_L Q$ if $Q=PS$ for some $S\in B_n^+$;
$P\le_R Q$ if $Q=SP$ for some $S\in B_n^+$.
\end{definition}

It is known that the posets $(B_n^+,\le_L)$ and $(B_n^+,\le_R)$ are lattices.
Let $\wedgeL$ and $\veeL$ (respectively, $\wedgeR$ and $\veeR$)
denote the gcd and the lcm with respect to $\le_L$ (respectively, $\le_R$).
For positive braids $P_1$ and $P_2$, the gcd $P_1\wedgeR P_2$ and
the lcm $P_1\veeR P_2$ are characterized
by the following properties:
\begin{itemize}
\item $P_1=Q_1(P_1\wedgeR P_2)$ and $P_2=Q_2(P_1\wedgeR P_2)$ for some
$Q_1,Q_2\in B_n^+$ with $Q_1\wedgeR Q_2=1$;
\item $P_1\veeR P_2=R_1P_1=R_2P_2$ for some
$R_1, R_2\in B_n^+$ with $R_1\wedgeL R_2=1$.
\end{itemize}
The partial orders $\le_L$ and $\le_R$, and thus the lattice
structures in $B_n^+$ can be extended to $B_n$ as follows:
for $\alpha,\beta\in B_n$,
$\alpha\le_L\beta$ if $\beta=\alpha P$ for some $P\in B_n^+$;
$\alpha\le_R\beta$ if $\beta=P\alpha$ for some $P\in B_n^+$.

\begin{definition}
The braid $\Delta=(\sigma_1\cdots\sigma_{n-1})(\sigma_1\cdots\sigma_{n-2})
\cdots(\sigma_1\sigma_2)\sigma_1$ is called the \emph{fundamental braid}
(or the \emph{Garside element}).
Let $\D=\{A\in B_n^+: A\le_L\Delta\}$.
The elements of $\D$ are called \emph{permutation braids}
(or \emph{simple elements}).
\end{definition}

The fundamental braid $\Delta$ has the following properties:
$A\le_L\Delta$ if and only if $A\le_R\Delta$ for $A\in B_n^+$;
$\Delta\le_L P$ if and only if $\Delta\le_R P$ for $P\in B_n^+$;
$\sigma_i\le_L \Delta$ and
$\sigma_i\Delta=\Delta\sigma_{n-i}$ for $i=1,\ldots,n-1$.
Permutation $n$-braids are in one-to-one correspondence
with $n$-permutations:
for an $n$-permutation $\theta$,
the diagram (in $[0,1]\times \mathbb R$) of the corresponding braid is obtained by
connecting
$(1,i)\in\{1\}\times \mathbb R$ to $(0,\theta(i))\in\{0\}\times\mathbb R$
by a straight line for each $i=1,\ldots,n$ and then
making the $i$-th strand lie
above the $j$-th strand whenever $i<j$.

For $P\in B_n^+$, let
$\LM(P)=P\wedge_L \Delta$ and $\RM(P)=P\wedge_R \Delta$.
It is known that for $P,Q\in B_n^+$,
$$\LM(PQ)=\LM(P\LM(Q))\quad\mbox{and}\quad \RM(PQ)=\RM(\RM(P)Q).$$

For $\alpha\in B_n$, there are integers $u\le v$ such that $\Delta^u\le_L
\alpha\le_L \Delta^v$. Let
$\inf(\alpha)=\max\{u\in\Z:\Delta^u\le_L\alpha\}$ and
$\sup(\alpha)=\min\{v\in\Z:\alpha\le_L\Delta^v\}$.

\begin{definition}
The expression $\Delta^u A_1\cdots A_m$
is called the \emph{left (respectively, right) normal form} of $\alpha$
if $u=\inf(\alpha)$, $A_i\in\D\setminus\{1,\Delta\}$ and $\LM(A_i\cdots A_m)=A_i$
(respectively, $\RM(A_1\cdots A_i)=A_i$) for $i=1,\ldots,m$.
\end{definition}

\begin{definition}
For $P\in B_n^+$, the \emph{starting set} $S(P)$ and
the \emph{finishing set} $F(P)$ of $P$ are defined as
$$S(P) = \{i\mid \sigma_i\le_L P\}\quad\mbox{and}\quad
F(P) = \{i\mid \sigma_i\le_R P\}.$$
\end{definition}

The following properties are well known~\cite{Thu92,EM94}.
\begin{lemma}\label{lemma:starting_set}
\begin{itemize}
\item[(i)]
For a positive braid $P$, $S(\LM(P))=S(P)$.
\item[(ii)]
If $A$ is a permutation braid with induced permutation $\theta$,
$$
S(A)=\{i\mid \theta^{-1}(i)>\theta^{-1}(i+1)\}\quad\mbox{and}\quad
F(A)=\{i\mid \theta(i)>\theta(i+1)\}.
$$
\item[(iii)]
For permutation braids $A$ and $B$, the expression
$AB$ is in left (respectively, right) normal form
if and only if $F(A)\supset S(B)$ (respectively, $F(A)\subset S(B)$).
\end{itemize}
\end{lemma}

By Thurston~\cite{Thu92},  an $n$-braid $\alpha$ has a unique expression
$$
\alpha=P^{-1}Q,
$$
where $P,Q\in B_n^+$ and $P\wedgeL Q=1$.
We call it the \emph{np}-form of $\alpha$.
Similarly, we define the \emph{pn}-form of $\alpha$
as $\alpha=PQ^{-1}$,
where $P,Q\in B_n^+$ and $P\wedgeR Q=1$.

Let $\tau$ be the inner automorphism of $B_n$
defined by $\tau(\sigma_i)=\sigma_{n-i}$.
Then $\Delta^{-1}\alpha\Delta=\tau(\alpha)$ for $\alpha\in B_n$.
The following is known~\cite[Lemma 2.3]{Cha95}.

\begin{lemma}\label{lemma:np-form}
Let $P, Q\in B_n^{+}$.
For $A\in\D$, let $\bar{A}=\Delta A^{-1}$.
\begin{enumerate}
\item[(i)]
Let $P=A_mA_{m-1}\cdots A_1$ and $Q=A_{m+1}A_{m+2}\cdots A_l$
be in left normal forms.
If $P^{-1}Q$ is in $np$-form, then
$\Delta^{-m} \tau^{1-m}(\bar A_1)\cdots \tau^{-1}(\bar A_{m-1}) \bar A_m
A_{m+1}\cdots A_l$
is the left normal form of $P^{-1}Q$.

\item[(ii)]
Let $P=A_1 A_2 \cdots A_m$ and $Q=A_l A_{l-1} \cdots A_{m+1}$
be in right normal forms.
If $PQ^{-1}$ is in $pn$-form, then
$\Delta^{m-l} \tau^{m-l}(A_1)\cdots\tau^{m-l}(A_m)\tau^{m-l+1}(\bar A_{m+1})
\cdots\tau^{-1}(\bar A_{l-1}) \bar A_{l}$
is the right normal form of $PQ^{-1}$.
\end{enumerate}
\end{lemma}

\subsection{Conjugacy problem in braid groups}
Let $\Delta^u A_1\cdots A_m$ be the left normal form
of $\alpha\in B_n$.
The \emph{cycling} $\c(\alpha)$ and the \emph{decycling} $\d(\alpha)$ are
defined by
\begin{eqnarray*}
\c(\alpha)&=&\Delta^u A_2\cdots A_m\tau^{-u}(A_1);\\
\d(\alpha)&=&\Delta^u \tau^u(A_m)A_1\cdots A_{m-1}.
\end{eqnarray*}

Let $[\alpha]$ denote the conjugacy class of $\alpha$.
Let $\infs(\alpha)=\max\{\inf(\beta):\beta\in[\alpha]\}$
and $\sups(\alpha)=\min\{\sup(\beta):\beta\in[\alpha]\}$.

\begin{definition}
For $\alpha\in B_n$, the \emph{super summit set} $[\alpha]^S$,
the \emph{ultra summit set} $[\alpha]^U$ and
the \emph{stable super summit set} $[\alpha]^\ST$
of $\alpha$ are defined as follows:
\begin{eqnarray*}
[\alpha]^S &=&\{\beta\in [\alpha]:
    \inf(\beta)=\infs(\alpha),\
    \sup(\beta)=\sups(\alpha)\};\\{}
[\alpha]^U &=&\{\beta\in [\alpha]^S:
    \mbox{$\c^m(\beta)=\beta$ for some $m\ge1$}\};\\{}
[\alpha]^\ST &=&\{\beta\in [\alpha]^S:
    \beta^m\in[\alpha^m]^S\
    \mbox{for all $m\ge1$}\}.
\end{eqnarray*}
\end{definition}

By definition, $[\alpha]^U$ and $[\alpha]^\ST$ are subsets of $[\alpha]^S$.

\begin{theorem}\label{thm:ConjAlgorithm}
Let $\alpha\in B_n$.
\begin{enumerate}
\item[(i)]
If\/ $\c^m(\alpha)=\alpha$ for some $m\ge 1$,
then $\inf(\alpha)=\infs(\alpha)$.

\item[(ii)]
If\/ $\d^m(\alpha)=\alpha$ for some $m\ge 1$,
then $\sup(\alpha)=\sups(\alpha)$.

\item[(iii)]
$\c^m\d^l(\alpha)\in[\alpha]^U$ for some $m,l\ge 0$.

\item[(iv)]
Both $[\alpha]^S$ and $[\alpha]^U$ are finite and nonempty.

\item[(v)] If\/ $\beta\in [\alpha]^S$,
then $\c(\beta),\d(\beta),\tau(\beta)\in [\alpha]^S$.
The same is true for $[\alpha]^U$.

\item[(vi)]
If\/ $\beta\in [\alpha]^S$, then
$\c(\d(\alpha))=\d(\c(\alpha))$,
$\tau(\c(\beta))=\c(\tau(\beta))$ and $\tau(\d(\beta))=\d(\tau(\beta))$.

\item[(vii)]
For $\beta,\beta'\in [\alpha]^S$,
there is a finite sequence
$$\beta=\beta_0\to \beta_1\to\cdots\to \beta_m=\beta'$$
such that for $i=0,\ldots,m-1$, $\beta_i\in [\alpha]^S$ and
$\beta_{i+1}=A_i\beta_iA_i^{-1}$ for some $A_i\in\D $.
The same is true for $[\alpha]^U$.
\end{enumerate}
\end{theorem}

For the results on stable super summit sets, see~\cite{LL06a,LL06b}.
For $\beta\in[\alpha]^S$, let
\begin{eqnarray*}
C^S(\beta) &=& \{P\in B_n^+:P^{-1}\beta P\in[\beta]^S\};\\
C^U(\beta) &=& \{P\in B_n^+:P^{-1}\beta P\in[\beta]^U\}.
\end{eqnarray*}
Both $C^S(\beta)$ and $C^U(\beta)$ are closed under $\wedgeL$
by Franco and Gonz\'alez-Meneses~\cite{FG03} and
Gebhardt~\cite{Geb05}, respectively.
The closedness under $\wedgeL$ makes the conjugacy algorithm more efficient.

For a nonempty subset $\mathcal V$ of $B_n^+$,
we call an element $P\in \mathcal V$
the \emph{$\le_R$-minimal element} of $\mathcal V$
if $P\le_R Q$ for all $Q\in \mathcal V$.
By definition, the $\le_R$-minimal element is unique if it exists.
If $\mathcal V$ is closed under $\wedge_R$, then $\mathcal V$
has the $\le_R$-minimal element.

\medskip
The following notions are useful in studying powers~\cite{LL07a,LL06b}.
For $\alpha\in B_n$, let
$$
\INF(\alpha)=\lim_{m\to\infty}\frac{\inf(\alpha^m)}m\quad\mbox{and}\quad
\SUP(\alpha)=\lim_{m\to\infty}\frac{\sup(\alpha^m)}m.
$$
The following lists important properties of $\INF(\cdot)$ and $\SUP(\cdot)$.
See Lemmas 3.2, 3.3, Theorem 3.13 in \cite{LL07a}, and
Corollary 3.5 in \cite{LL06b}.
\begin{theorem}\label{thm:INF}
Let $\alpha\in B_n$.
\begin{enumerate}
\item[(i)]
$\INF(\gamma\alpha\gamma^{-1})=\INF(\alpha)$ and
$\SUP(\gamma\alpha\gamma^{-1})=\SUP(\alpha)$ for all $\gamma\in B_n$.

\item[(ii)]
$\INF(\alpha^m)=m\INF(\alpha)$ and
$\SUP(\alpha^m)=m\SUP(\alpha)$ for all $m\ge 1$.

\item[(iii)]
$\infs(\alpha)\le\INF(\alpha)<\infs(\alpha)+1$ and
$\sups(\alpha)-1<\SUP(\alpha)\le\sups(\alpha)$.

\item[(iv)]
$\INF(\alpha)$ and $\SUP(\alpha)$ are rational of the form $p/q$ for some
integers $p,q$ with $1\le q\le n(n-1)/2$.

\end{enumerate}
\end{theorem}

\subsection{Duality between cycling and decycling}
In many aspects, the cycling and the decycling are dual to each other.
We define a variant of the cycling as follows
so that the duality is more clear. See Lemmas~\ref{thm:duality}
and \ref{thm:duality2}.
\begin{definition}
For $\alpha\in B_n$, define $\c_0(\alpha)=\tau^{-1}(\c(\alpha))$.
\end{definition}

Since $\tau^2(\beta)=\beta$ and $\tau(\c(\beta))=\c(\tau(\beta))$
for $\beta\in[\alpha]^S$,
we can replace $\c$ with $\c_0$ in Theorem~\ref{thm:ConjAlgorithm}
and in the definition of $[\alpha]^U$.
In particular, for an element $\beta\in[\alpha]^S$,
$\beta$ belongs to the ultra summit set $[\alpha]^U$ if and only if
$\c_0^m(\beta)=\beta$ for some $m\ge 1$.

\begin{lemma}~\label{thm:duality}
Let $\Delta^u A_1\cdots A_m$ be the left normal form of $\alpha\in B_n$.
\begin{enumerate}
\item[(i)]
The set $\{P\in B_n^+: \inf(P\alpha)>\inf(\alpha)\}$ is
nonempty and closed under $\wedgeR$.
The $\le_R$-minimal element $A$ of this set is the permutation braid
$\tau^{-u}(\Delta A_1^{-1})$ and
satisfies $\c_0(\alpha)= A\alpha A^{-1}$.

\item[(ii)]
The set $\{P\in B_n^+: \sup(\alpha P^{-1})<\sup(\alpha)\}$ is
nonempty and closed under $\wedgeR$.
The $\le_R$-minimal element $A$ of this set is the permutation braid $A_m$
and satisfies $\d(\alpha)=A\alpha A^{-1}$.
\end{enumerate}
\end{lemma}

\begin{proof}
We prove only (i) since (ii) can be proved similarly.
Nonemptiness of $\{P\in B_n^+: \inf(P\alpha)>\inf(\alpha)\}$ is clear.
Note that
\begin{itemize}
\item
$(\beta\alpha)\wedge_R(\gamma\alpha)=(\beta\wedge_R \gamma)\alpha$
for all $\alpha,\beta,\gamma\in B_n$;
\item
$\inf(\alpha\wedge_R \beta)=\min\{\inf(\alpha),\inf(\beta)\}$
for all $\alpha,\beta\in B_n$.
\end{itemize}

If $\inf(P\alpha)>\inf(\alpha)$ and $\inf(Q\alpha)>\inf(\alpha)$ for
positive braids $P$ and $Q$, then
$$
\inf((P\wedge_R Q)\alpha)
=\inf((P\alpha)\wedge_R(Q\alpha))
=\min\{\inf(P\alpha),\inf(Q\alpha)\}
>\inf(\alpha).
$$
Therefore, the set $\{P\in B_n^+: \inf(P\alpha)>\inf(\alpha)\}$ is
closed under $\wedge_R$.

It is easy to see that the $\le_R$-minimal element $A$ is
$\tau^{-u}(\Delta A_1^{-1})$ and, hence,
\begin{eqnarray*}
A\alpha A^{-1} &=& (\Delta \tau^{-u}(A_1^{-1}))\
(\Delta^u A_1\cdots A_m)\
(\tau^{-u}(A_1)\Delta^{-1})\\
&=&\Delta\
 (\Delta^u A_2\cdots A_m \tau^{-u}(A_1))\
 \Delta^{-1}
=\Delta\c(\alpha)\Delta^{-1}
=\tau^{-1}(\c(\alpha))\\
&=&\c_0(\alpha).
\end{eqnarray*}
\vskip-\baselinestretch\baselineskip
\end{proof}

\begin{definition}
For $\alpha\in B_n$, the set
$$
[\alpha]^U_\d=\{\beta\in[\alpha]^S:\d^m(\beta)=\beta
\mbox{ for some $m\ge 1$}\}
$$
is called the \emph{ultra summit set
of $\alpha$ with respect to decycling}.
\end{definition}

The following lemma is easy to prove, so we omit the proof.
It shows that there is a duality between
$\c_0(\cdot)\leftrightarrow \d(\cdot)$,
$\inf(\cdot)\leftrightarrow \sup(\cdot)$
and $[\,\cdot\,]^U\leftrightarrow[\,\cdot\,]^U_\d$.

\begin{lemma}~\label{thm:duality2}
Let $\alpha\in B_n$.
\begin{enumerate}
\item[(i)]
$\inf(\alpha)=-\sup(\alpha^{-1})$ and\/ $\infs(\alpha)=-\sups(\alpha^{-1})$.

\item[(ii)]
$\c_0(\alpha)=(\d(\alpha^{-1}))^{-1}$.

\item[(iii)]
$\beta\in[\alpha]^S$ if and only if\/ $\beta^{-1}\in[\alpha^{-1}]^S$.

\item[(iv)]
$\beta\in[\alpha]^U$ if and only if\/ $\beta^{-1}\in[\alpha^{-1}]^U_\d$.
\end{enumerate}
\end{lemma}

\section{Braids sending a standard curve to a standard curve}
\label{sec:standard}

In this section we study the normal form of braids that send
a standard curve system to a standard curve system.
We collect basic properties of such braids in Lemma~\ref{thm:decom},
from which the other results of this section follow easily.

We start by defining some notions.
Throughout the paper, we do not distinguish the curves
and the isotopy classes of curves.

\begin{figure}
\includegraphics[scale=.6]{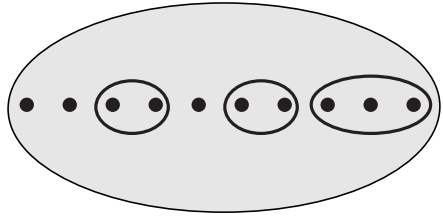}
\vskip -2mm
\caption{The unnested standard curve system $\C_\n$ for  $\n=(1,1,2,1,2,3)$}
\label{fig:unnested}
\end{figure}

\begin{definition}\label{def:curves}
A curve system means a finite collection of disjoint simple closed curves.
A simple closed curve in $D_n$ is said to be \emph{essential}\/
if it is homotopic neither to a point nor to a puncture nor to
the boundary.
An essential curve system in $D_n$ is said to be \emph{standard}\/
if each component is isotopic to a round circle
centered at the real axis as in Figure~\ref{fig:standard}.
It is said to be \emph{unnested} if none of its components
encloses another component.
See Figure~\ref{fig:unnested}.
\end{definition}

\begin{definition}
The $n$-braid group $B_n$ acts on the set of curve systems in $D_n$.
Let $\alpha*\C$ denote the left action of $\alpha\in B_n$
on the curve system $\C$ in $D_n$.
An $n$-braid $\alpha$ is said to be \emph{reducible}
if $\alpha*\C=\C$ for some essential curve system $\C$ in $D_n$.
Such a curve system $\C$ is called a \emph{reduction system} of $\alpha$.
\end{definition}

The unnested standard curve systems in $D_n$ are
in one-to-one correspondence
with the $k$-compositions of $n$ for $2\le k\le n-1$.
Recall that an ordered $k$-tuple $\n=(n_1,\ldots,n_k)$
is a $k$-composition of $n$ if
$n_i\ge 1$ for each $i$ and $n=n_1+\cdots+n_k$.

\begin{definition}
For a composition $\n=(n_1,\ldots,n_k)$ of $n$,
let $\C_\n$ denote the curve system $\cup_{n_i\ge 2}C_i$,
where $C_i$ is the standard curve enclosing
$\{m: \sum_{j=1}^{i-1}n_j< m\le \sum_{j=1}^{i}n_j\}$.
See Figure~\ref{fig:unnested}.
\end{definition}

The $k$-braid group $B_k$ acts on the set of $k$-compositions
of $n$ via the induced permutations: for a $k$-composition
$\n=(n_1,\cdots,n_k)$ and
$\alpha_0\in B_k$ with induced permutation $\theta$,
$\alpha_0*\n=(n_{\theta^{-1}(1)},\ldots,n_{\theta^{-1}(k)})$.

\begin{definition}
Let $\n=(n_1,\cdots,n_k)$ be a composition of $n$.
\begin{itemize}
\item
Let $\alpha_0=l_1\cup\cdots\cup l_k$ be a $k$-braid
with $l_i\cap(\{1\}\times D^2)=\{(1,i)\}$ for each $i$.
Note that the strands of $\alpha_0$ are numbered from bottom to top
at its right end.
We define $\myangle{\alpha_0}_\n$ as the $n$-braid obtained
from $\alpha_0$ by taking $n_i$ parallel copies of $l_i$ for each $i$.
\item
Let $\alpha_i\in B_{n_i}$ for $i=1,\ldots,k$.
We define $(\alpha_1\oplus\cdots\oplus\alpha_k)$ as the $n$-braid
$\alpha_1'\alpha_2'\cdots\alpha_k'$, where $\alpha_i'$
is the image of $\alpha_i$ under the homomorphism
$B_{n_i}\to B_n$ defined by $\sigma_j\mapsto\sigma_{n_1+\cdots+n_{i-1}+j}$.
\end{itemize}
\end{definition}

\begin{figure}
\begin{tabular}{ccccc}
\includegraphics[scale=1]{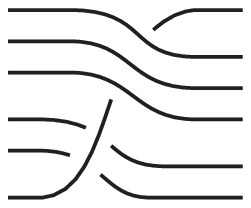} &\qquad&
\includegraphics[scale=1]{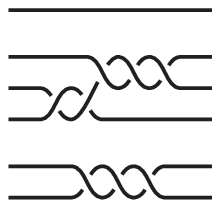} &\qquad&
\includegraphics[scale=1]{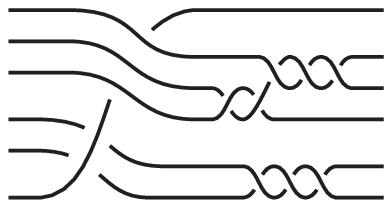} \\
\small (a) $\myangle{\sigma_1^{-1}\sigma_2}_\n$&&
\small (b) $(\sigma_1^3\oplus\sigma_1^{-2}\sigma_2^3\oplus 1)$&&
\small (c) $\myangle{\sigma_1^{-1}
\sigma_2}_\n(\sigma_1^3\oplus\sigma_1^{-2}\sigma_2^3\oplus 1)$
\end{tabular}
\caption{$\n=(2,3,1)$}\label{fig:copy}
\end{figure}

We will use the notation
$\alpha=\myangle{\alpha_0}_\n(\alpha_1\oplus\cdots\oplus\alpha_k)$
throughout the paper. See Figure~\ref{fig:copy}.

\begin{lemma}\label{thm:decom}
Let\/ $\n=(n_1,\ldots,n_k)$ be a composition of\/ $n$.
\begin{enumerate}
\item[(i)]
The expression $\alpha=\myangle{\alpha_0}_\n(\alpha_1\oplus\cdots\oplus\alpha_k)$
is unique, i.e., if\/
$\myangle{\alpha_0}_\n(\alpha_1\oplus\cdots\oplus\alpha_k)
=\myangle{\beta_0}_\n(\beta_1\oplus\cdots\oplus\beta_k)$,
then $\alpha_i=\beta_i$ for $i=0,\ldots,k$.

\item[(ii)]
If $\alpha=\myangle{\alpha_0}_\n(\alpha_1\oplus\cdots\oplus\alpha_k)$,
then $\alpha*\C_\n$ is standard and, further, $\alpha*\C_\n=\C_{\alpha_0*\n}$.
Conversely, if\/ $\alpha*\C_\n$ is standard, then $\alpha$ can be expressed as
$\alpha=\myangle{\alpha_0}_\n(\alpha_1\oplus\cdots\oplus\alpha_k)$.

\item[(iii)]
Let $\alpha=\myangle{\alpha_0}_\n(\alpha_1\oplus\cdots\oplus\alpha_k)$.
If all $\alpha_i$'s are positive (respectively, permutation  and fundamental)
braids, then so is $\alpha$.

\item[(iv)]
$\myangle{\alpha_0}_\n(\alpha_1\oplus\cdots\oplus\alpha_k)
= (\alpha_{\theta^{-1}(1)}\oplus\cdots\oplus\alpha_{\theta^{-1}(k)})
\myangle{\alpha_0}_\n$,
where $\theta$ is the induced permutation of\/ $\alpha_0$.

\item[(v)]
$\myangle{\alpha_0\beta_0}_\n
=\myangle{\alpha_0}_{\beta_0*\n}\myangle{\beta_0}_\n$.

\item[(vi)]
$(\myangle{\alpha_0}_\n)^{-1}=\myangle{\alpha_0^{-1}}_{\alpha_0*\n}$.

\item[(vii)]
$
(\alpha_1\beta_1\oplus\cdots\oplus\alpha_k\beta_k)
=(\alpha_1\oplus\cdots\oplus\alpha_k)
(\beta_1\oplus\cdots\oplus\beta_k)
$

\item[(viii)]
$(\alpha_1\oplus\cdots\oplus\alpha_k)^{-1}
=(\alpha_1^{-1}\oplus\cdots\oplus\alpha_k^{-1})$.

\item[(ix)]
Let $A_0$ and $B_0$ be permutation $k$-braids.
$A_0B_0$ is in left (respectively, right) normal form
if and only if\/ $\myangle{A_0}_{B_0*\n}\,\myangle{B_0}_\n$
is in left (respectively, right) normal form.

\item[(x)]
Let $P_i$, $i=0,\ldots,k$, be positive braids with appropriate braid indices.
Let $A_i=\LM(P_i)$ and $B_i=\RM(P_i)$ for $i=0,\ldots,k$. Then
\begin{eqnarray*}
\LM((P_1\oplus\cdots\oplus P_k)\myangle{P_0}_\n)
&=&(A_1\oplus\cdots\oplus A_k)\myangle{A_0}_{(A_0^{-1}P_0)*\n};\\
\RM(\myangle{P_0}_\n(P_1\oplus\cdots\oplus P_k))
&=&\myangle{B_0}_\n(B_1\oplus\cdots\oplus B_k).
\end{eqnarray*}
\end{enumerate}
\end{lemma}

\begin{proof}
The statements from (i) to (viii) are easy to prove.
Let us prove (ix) and (x).

\smallskip
(ix) \
Let $B_0*\n=(n_1',\ldots,n_k')$ and $N_i=n_1'+\cdots+n_i'$
for $i=1,\ldots,k$. Then,
\begin{eqnarray*}
F(\myangle{A_0}_{B_0*\n})
&=&\{N_i: i\in F(A_0)\};\\
S(\myangle{B_0}_\n)
&=&\{N_i: i\in S(B_0)\}.
\end{eqnarray*}
Hence, $F(A_0)\supset S(B_0)$ if and only if
$F(\myangle{A_0}_{B_0*\n})\supset S(\myangle{B_0}_\n)$,
and $F(A_0)\subset S(B_0)$ if and only if
$F(\myangle{A_0}_{B_0*\n})\subset S(\myangle{B_0}_\n)$.

\smallskip
(x) \
We prove only the second identity.
The first one can be proved in a similar way.
It is easy to see that
$\RM(\myangle{P_0}_\n)=\myangle{B_0}_\n$ by (ix)
and that $\RM(P_1\oplus\cdots\oplus P_k)=(B_1\oplus\cdots\oplus B_k)$.
Let $\theta$ be the induced permutation of $B_0$.
Then, by (iv)
\begin{eqnarray*}
\lefteqn{\RM(\myangle{P_0}_\n (P_1\oplus\cdots\oplus P_k))
  = \RM(\RM(\myangle{P_0}_\n) (P_1\oplus\cdots\oplus P_k))}\\
&=& \RM(\myangle{B_0}_\n (P_1\oplus\cdots\oplus P_k))
  = \RM((P_{\theta^{-1}(1)}\oplus\cdots\oplus P_{\theta^{-1}(k)})\myangle{B_0}_\n)\\
&=& \RM(\RM(P_{\theta^{-1}(1)}\oplus\cdots\oplus P_{\theta^{-1}(k)})\myangle{B_0}_\n)
  = \RM((B_{\theta^{-1}(1)}\oplus\cdots\oplus B_{\theta^{-1}(k)})\myangle{B_0}_\n)\\
&=& \RM(\myangle{B_0}_\n(B_1\oplus\cdots\oplus B_k))
  = \myangle{B_0}_\n(B_1\oplus\cdots\oplus B_k).
\end{eqnarray*}
The last equality holds since $\myangle{B_0}_\n(B_1\oplus\cdots\oplus B_k)$
is a permutation braid by (iii).
\end{proof}

Let $\brindex(\alpha)$ denote the braid index of $\alpha$.

\begin{lemma}\label{thm:infofredbr}
Let $\alpha=\myangle{\alpha_0}_\n(\alpha_1\oplus\cdots
\oplus\alpha_k)\in B_n$.
\begin{enumerate}
\item[(i)]
$\inf(\alpha)=\min\{\inf(\alpha_i):
i=0,\ldots,k,\ \brindex(\alpha_i)\ge 2\}$.

\item[(ii)]
$\sup(\alpha)=\max\{\sup(\alpha_i):
i=0,\ldots,k,\ \brindex(\alpha_i)\ge 2\}$.

\item[(iii)]
$\alpha$ is a positive (respectively, permutation  and fundamental)
braid if and only if each $\alpha_i$ is a
positive (respectively, permutation  and fundamental)
braid for $i=0,\ldots,k$.
\end{enumerate}
\end{lemma}

\begin{proof}
(i) \
Let $r=\min\{\inf(\alpha_i):i=0,\ldots,k,\ \brindex(\alpha_i)\ge 2\}$.
Set $n_0=k$.
For $i=0,\ldots,k$, let $\alpha_i=\Delta_i^rP_i$,
where $\Delta_i$ is the fundamental braid of $B_{n_i}$ and $P_i\in B_{n_i}^+$.
Let $P=\myangle{P_0}_\n(P_1\oplus\cdots\oplus P_k)$.
By Lemma~\ref{thm:decom}~(iv), (v) and (vii),
\begin{eqnarray*}
\alpha
&=&\myangle{\Delta_0^rP_0}_\n(\Delta_1^rP_1\oplus\cdots\oplus \Delta_k^rP_k)\\
&=&\myangle{\Delta_0^r}_{P_0*\n}
  \myangle{P_0}_\n
  (\Delta_1^r\oplus\cdots\oplus \Delta_k^r)
  (P_1\oplus\cdots\oplus P_k)\\
&=&\myangle{\Delta_0^r}_{P_0*\n}
  (\Delta_{\theta^{-1}(1)}^r\oplus\cdots\oplus \Delta_{\theta^{-1}(k)}^r)
  \myangle{P_0}_\n
  (P_1\oplus\cdots\oplus P_k)
\end{eqnarray*}
where $\theta$ is the induced permutation of $P_0$.
Since $P_0*\n=(n_{\theta^{-1}(1)},\ldots,n_{\theta^{-1}(k)})$,
we have
$\myangle{\Delta_0^r}_{P_0*\n}
(\Delta_{\theta^{-1}(1)}^r\oplus\cdots\oplus \Delta_{\theta^{-1}(k)}^r)
=\Delta^r$,
and hence $\alpha=\Delta^rP$.
Since $\inf(P_i)=0$ for some $P_i$ with $\brindex(P_i)\ge 2$,
$\RM(P)\ne\Delta$ by Lemma~\ref{thm:decom} (x).
Therefore $\inf(\alpha)=r$.

\smallskip
(ii) \
Since $\sup(\alpha)=-\inf(\alpha^{-1})$ by Lemma~\ref{thm:duality2}~(i) and
$\alpha^{-1}=(\alpha_1^{-1}\oplus\cdots\oplus\alpha_k^{-1})
\myangle{\alpha_0^{-1}}_{\alpha_0*\n}$ by Lemma~\ref{thm:decom} (vi) and (viii),
the assertion follows from (i).

\smallskip
(iii) \
Note that a braid $\beta$ is a
positive (respectively, permutation  and fundamental) braid
if and only if $\inf(\beta)\ge 0$
(respectively, $0\le \inf(\beta)\le\sup(\beta)\le 1$
and $\inf(\beta)=\sup(\beta)=1$).
Therefore, the assertion follows from (i) and (ii) and Lemma~\ref{thm:decom}~(iii).
\end{proof}

\begin{lemma}\label{thm:normalNstandardLemma}
Let $\C$ be a standard curve system in $D_n$ and $P\in B_n^+$
such that $P*\C$ is standard.
\begin{enumerate}
\item[(i)]
If $P=QA$ and $A=\RM(P)$, then $A*\C$ is standard.

\item[(ii)]
If $P=AQ$ and $A=\LM(P)$, then $Q*\C$ is standard.
\end{enumerate}
\end{lemma}

\begin{proof}
A curve system is standard if and only if each of its components is standard.
Hence, we may assume that the given standard curve system $\C$ is unnested.
Let $\C=\C_\n$ for a composition $\n=(n_1,\ldots,n_k)$ of $n$.

\smallskip
(i)\
$P=\myangle{P_0}_\n(P_1\oplus\cdots \oplus P_k)$
for some positive braids $P_i$, $i=0,\ldots,k$,
by Lemmas~\ref{thm:decom}~(ii) and~\ref{thm:infofredbr}~(iii).
By Lemma~\ref{thm:decom} (x),
$A=\RM(P)=\myangle{\RM(P_0)}_\n(\RM(P_1)\oplus\cdots\oplus \RM(P_k))$.
By Lemma~\ref{thm:decom} (ii), $A*\C$ is standard.

\smallskip
(ii)
$P=(P_1\oplus\cdots \oplus P_k)\myangle{P_0}_\n$
for some positive braids $P_i$, $i=0,\ldots,k$,
by Lemmas~\ref{thm:decom}~(ii), (iv) and~\ref{thm:infofredbr}~(iii).
Let $A_i=\LM(P_i)$ for $i=0,\ldots,k$.
Then $A=\LM(P)=(A_1\oplus\cdots\oplus A_k)\myangle{A_0}_{(A_0^{-1}P_0)*\n}$
by Lemma~\ref{thm:decom} (x).
By Lemma~\ref{thm:decom} (vi) and (viii),
$A^{-1}=\myangle{A_0^{-1}}_{P_0*\n}
(A_1^{-1}\oplus\cdots\oplus A_k^{-1})$.
By Lemma~\ref{thm:decom} (ii),
\begin{eqnarray*}
Q*\C_\n&=&(A^{-1}P)*\C_\n=A^{-1}*(P*\C_\n)=A^{-1}*\C_{P_0*\n}\\
&=& \bigl(\myangle{A_0^{-1}}_{P_0*\n}
(A_1^{-1}\oplus\cdots\oplus A_k^{-1})\bigr)
*\C_{P_0*\n}
= \C_{(A_0^{-1}P_0)*\n}.
\end{eqnarray*}
Hence $Q*\C$ is standard.
\end{proof}

\begin{theorem}\label{thm:normalNstandard}
Let $\C$ be a standard curve system in $D_n$
and $\Delta^u A_1\cdots A_m$ be
the (left or right) normal form of\/ $\alpha\in B_n$.
If\/ $\alpha*\C$ is standard,
then so is $(A_i\cdots A_m)*\C$ for $i=1,\ldots,m$.
\end{theorem}

\begin{proof}
It is an immediate consequence of Lemma~\ref{thm:normalNstandardLemma},
because $(A_1\cdots A_m)*\C=\Delta^{-u}*(\alpha*\C)$ is standard.
\end{proof}

\begin{figure}
\raisebox{10mm}{$\alpha=$\ \ }\includegraphics[scale=1]{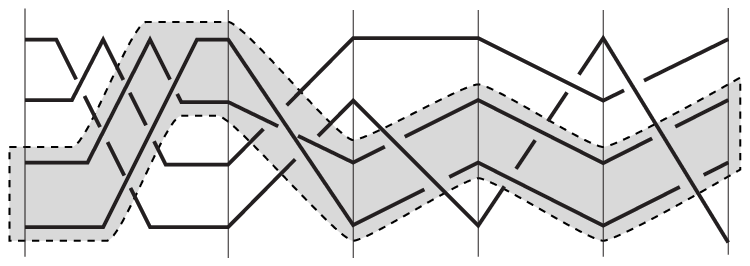}
\vskip -3mm
\caption{The 4-braid $\alpha$, whose normal form is
of the form $\Delta^{-1}A_1A_2A_3A_4$, sends
the standard curve system $\C_{(1,2,1)}$
to the standard curve system $\C_{(2,1,1)}$ as follows:
$
\C_{(2,1,1)}\stackrel{\,\,\Delta^{-1}}{\longleftarrow}
\C_{(1,1,2)}\stackrel{\,\,A_1}{\longleftarrow}
\C_{(2,1,1)}\stackrel{\,\,A_2}{\longleftarrow}
\C_{(1,2,1)}\stackrel{\,\,A_3}{\longleftarrow}
\C_{(2,1,1)}\stackrel{\,\,A_4}{\longleftarrow}
\C_{(1,2,1)}.
$
}
\label{fig:normalred}
\end{figure}

Roughly speaking, Theorem~\ref{thm:normalNstandard} says that
if a braid $\alpha$ sends a standard curve system to a standard curve system,
then so does each permutation braid in the normal form of $\alpha$
as in Figure~\ref{fig:normalred}.

\begin{corollary}\label{thm:BNG0}
Let $\Delta^u A_1\cdots A_m$ be the left normal form of\/
an $n$-braid $\alpha$.
If\/ $\alpha$ has a standard reduction system $\C$,
then $\c_0(\alpha)$, $\d(\alpha)$ and $\tau(\alpha)$
have standard reduction systems
$\tau^{-u}(\Delta A_1^{-1})*\C$, $A_m*\C$ and $\Delta^{-1}*\C$,
respectively.
\end{corollary}

\begin{proof}
$A_m*\C$ is standard by Theorem~\ref{thm:normalNstandard}.
By Lemma~\ref{thm:duality},
$$
\d(\alpha)*(A_m*\C)=(A_m\alpha A_m^{-1})*(A_m*\C)=A_m*(\alpha*\C)=A_m*\C.
$$
Therefore $\d(\alpha)$ has a standard reduction system $A_m*\C$.
In the same way,
$\tau(\alpha)$ and $\c_0(\alpha)$ have standard reduction systems
$\Delta^{-1}*\C$ and $\tau^{-u}(\Delta A_1^{-1})*\C$, respectively.
\end{proof}

\begin{corollary}\label{thm:BNG}
Let $\alpha$ be a reducible $n$-braid with a reduction system $\C$.
There exists an element $\beta$ of the ultra summit set $[\alpha]^U$ which has
a standard reduction system. Precisely,
there exists a positive braid $P$ such that
$\beta=P\alpha P^{-1}$ belongs to $[\alpha]^U$ and $P*\C$ is a standard
reduction system of\/ $\beta$.
\end{corollary}

\begin{proof}
Let $P_1$ be a positive $n$-braid such that $P_1*\C$ is standard.
Then $P_1\alpha P_1^{-1}$ has
the standard reduction system $P_1*\C$.
Take $l,m\ge 0$ such that $\beta=\c_0^l\d^m(P_1\alpha P_1^{-1})$
belongs to $[\alpha]^U$.
Lemma~\ref{thm:duality} and Corollary~\ref{thm:BNG0} say that
if $\gamma\in B_n$ has a standard reduction system $\C'$,
then there are permutation braids $A_1$ and $A_2$ such that
$\c_0(\gamma)=A_1\gamma A_1^{-1}$ and $\d(\gamma)=A_2\gamma A_2^{-1}$
have standard reduction systems $A_1 * \C'$ and $A_2 * \C'$, respectively.
Hence, we can find a positive $n$-braid $P_2$ such that
$\beta=P_2(P_1\alpha P_1^{-1})P_2^{-1}$ and
$P_2*(P_1*\C) = (P_2P_1)* \C$ is standard.
Let $P_2P_1=P$.
Then, $\beta=P\alpha P^{-1}$ and $\beta$ has the standard reduction system
$(P_2P_1)*\C=P*\C$.
\end{proof}

\begin{corollary}\label{thm:np-form}
Let\/ $\C$ be a standard curve system in $D_n$,
and let $\alpha*\C$ be standard for an $n$-braid $\alpha$.
\begin{enumerate}
\item[(i)]
If\/ $P^{-1}Q$ is the np-form of\/ $\alpha$, then $Q*\C$ is standard.

\item[(ii)]
If\/ $PQ^{-1}$ is the pn-form of\/ $\alpha$, then $Q^{-1}*\C$ is standard.
\end{enumerate}
\end{corollary}

\begin{proof}
By Lemma~\ref{lemma:np-form} and Theorem~\ref{thm:normalNstandard},
$Q*\C$ and $Q^{-1}*\C$ are standard.
\end{proof}

We remark that Theorem~\ref{thm:normalNstandard} and Corollary~\ref{thm:BNG}
were obtained also by Benardete, Guti\'errez and Nitecki~\cite[Theorems 5.7 and 5.8]{BGN95},
and that these two are enough to solve the reducibility problem
because there is an efficient algorithm that decides
whether a given braid has a standard reduction system or not
and finds one if it has~\cite{BGN93}.
However, Corollary~\ref{thm:BNG} guarantees only the
\emph{existence} of an element (in the ultra summit set of a reducible braid)
that has a standard reduction system.
To solve the reducibility problem using only Corollary~\ref{thm:BNG},
we have to compute all the elements in the ultra summit set.

\section{Standardizers of curve systems}

\begin{definition}
For an essential curve system $\C$ in $D_n$,
we define the \emph{standardizer} of $\C$ as the set
$$ \St(\C)= \{P\in B_n^+:\mbox{$P*\C$ is standard}\}. $$
\end{definition}

This section is devoted to the study of properties of standardizers.
Clearly, $\St(\C)$ is nonempty for any essential curve system $\C$.
Theorem~\ref{thm:closed} shows that standardizers are
sublattices of $B_n^+$,
hence they have unique $\le_R$-minimal elements.
The main result of this section is Theorem~\ref{thm:StMain}
that for any reduction system $\C$ of a reducible braid $\alpha$,
conjugating $\alpha$ by the $\le_R$-minimal element of $\St(\C)$
preserves the membership of the super summit set, ultra summit set
and stable super summit set.
Proposition~\ref{thm:FixStCur} and Corollary~\ref{cor:StCur}
show that the $\le_R$-minimal element of $\St(\C)$ does not entangle
any standard curve disjoint from $\C$.
Proposition~\ref{thm:MinElt} is a characterization
of the $\le_R$-minimal element of $\St(\C)$
in terms of normal form and lattice operations.

\begin{theorem}\label{thm:closed}
For an essential curve system $\C$ in $D_n$,
its standardizer $\St(\C)$ is closed under $\wedgeR$ and $\veeR$,
and hence a sublattice of $B_n^+$.
Therefore $\St(\C)$ contains a unique $\le_R$-minimal element.
\end{theorem}

\begin{proof} (See Figure~\ref{fig:lattice}.)
\begin{figure}
$
\xymatrix{
& P_1*\C \ar[ddl]_{R_1}\\ \mbox{}\\
(P_1\veeR P_2)*\C
&& (P_1\wedgeR P_2)*\C \ar[uul]_{Q_1}\ar[ddl]_{Q_2}
&& \C
\ar@/-1em/[ll]_{P_1\wedgeR P_2}
\ar@/_1em/[uulll]_{P_1}
\ar@/^1em/[ddlll]^{P_2}
\\ \mbox{}\\
& P_2*\C \ar[uul]_{R_2}
}
$
\vskip -2mm
\caption{$\C_1\stackrel{P}\longrightarrow\C_2$
means that $\C_2=P*\C_1$.}
\label{fig:lattice}
\end{figure}
Let $P_1,P_2\in\St(\C)$.
Let $P_1=Q_1(P_1\wedgeR P_2)$ and $P_2=Q_2(P_1\wedgeR P_2)$
for $Q_1,Q_2\in B_n^+$ with $Q_1\wedgeR Q_2=1$.
Then $P_2=Q_2(P_1\wedgeR P_2)=Q_2Q_1^{-1}P_1$, and
$Q_2Q_1^{-1}$ is in \emph{pn}-form.
Since $P_1*\C$ and $P_2*\C$ are standard and
$$P_2*\C=(Q_2Q_1^{-1})*(P_1*\C),$$
$Q_1^{-1}*(P_1*\C)=(P_1\wedgeR P_2)*\C$ is standard
by Corollary~\ref{thm:np-form} (ii).

Let $P_1\veeR P_2=R_1P_1=R_2P_2$ for $R_1,R_2\in B_n^+$ with $R_1\wedgeL R_2=1$.
Then $R_2^{-1}R_1P_1=P_2$, and $R_2^{-1}R_1$ is the \emph{np}-form.
Since $P_1*\C$ and $P_2*\C$ are standard and
$$P_2*\C=(R_2^{-1}R_1)*(P_1*\C),$$
$R_1*(P_1*\C)=(P_1\veeR P_2)*\C$ is standard
by Corollary~\ref{thm:np-form} (i).
\end{proof}

Let $\C$, $\C_1$ and $\C_2$ be essential curve systems
such that $\C=\C_1\cup\C_2$.
Then $\St(\C)\subset\St(\C_i)$ for $i=1,2$.
Let $P$, $P_1$ and $P_2$ be the $\le_R$-minimal elements of
$\St(\C)$, $\St(\C_1)$ and $\St(\C_2)$, respectively.
By Theorem~\ref{thm:closed}, $P_1\le_R P$ and $P_2\le_R P$,
hence $(P_1\veeR P_2)\le_R P$.
One may expect that $P=P_1\vee_R P_2$.
However, the following example shows that it is not true in general.

\begin{figure}
$
\xymatrix{
&
\includegraphics[scale=.6]{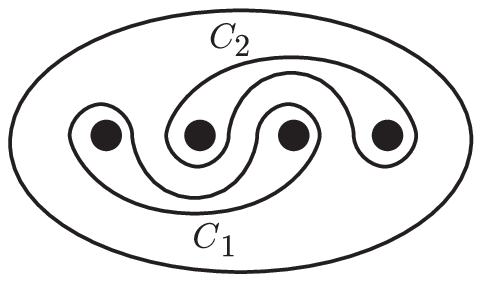}
\ar[dl]_{\sigma_1}
\ar[dr]^{\sigma_3}\\
\includegraphics[scale=.6]{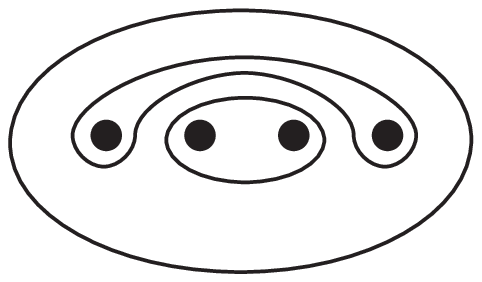}
\ar[dr]_{\sigma_2\sigma_3}
&&
\includegraphics[scale=.6]{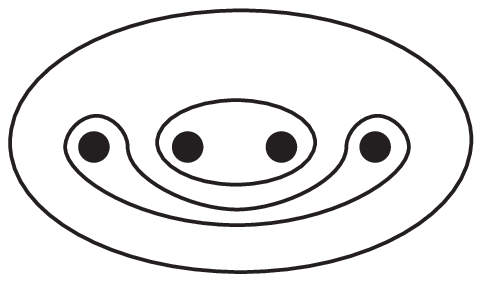}
\ar[dl]_{\sigma_2\sigma_1}\\
&\includegraphics[scale=.6]{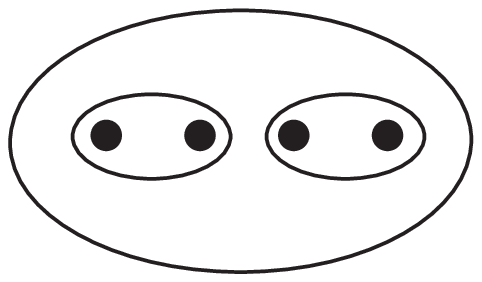}
}
$
\vskip -2mm
\caption{Standardization of a curve system}\label{fig:nolcm}
\end{figure}

\begin{example}
Let $C_1$ and $C_2$ be the curves in $D_4$ as in Figure~\ref{fig:nolcm}.
The $\le_R$-minimal elements of $\St(C_1)$, $\St(C_2)$ and $\St(C_1\cup C_2)$
are $\sigma_1$, $\sigma_3$ and $\sigma_2\sigma_1\sigma_3$, respectively.
Note that $\sigma_2\sigma_1\sigma_3$ is not equal to
$\sigma_1\veeR\sigma_3=\sigma_1\sigma_3$.
\end{example}

The following proposition shows that, when an essential curve $C$ in $D_n$
is standardized by the action of the $\le_R$-minimal element of $\St(C)$,
any other standard curve disjoint from $C$ remains standard.

\begin{proposition}\label{thm:FixStCur}
Let $C$ be an essential simple closed curve in $D_n$
and let $P$ be the $\le_R$-minimal element of $\St(C)$.
For any standard curve $C'$ in $D_n$ with $C\cap C'=\emptyset$,
the curve $P*C'$ is standard.
\end{proposition}

\begin{figure}
\tabcolsep=15pt
\begin{tabular}{c}
\includegraphics[scale=1.2]{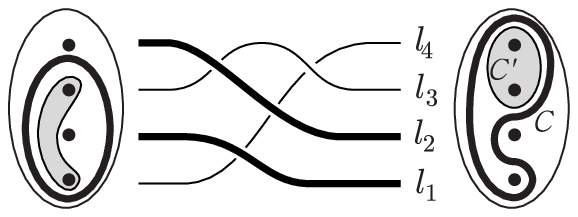}\\
(a) The braid diagram $K$ of $P$\\[6pt]
\includegraphics[scale=1.2]{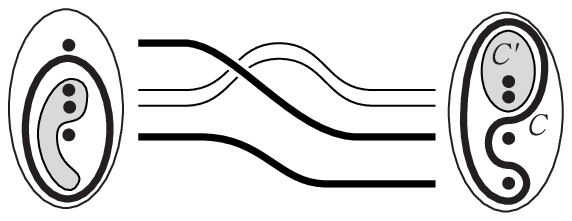}\quad
\raisebox{22pt}{\LARGE $=$}\quad
\includegraphics[scale=1.2]{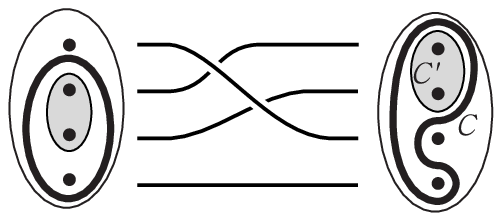}\\
(b) The braid diagram $L$ of $Q$
\end{tabular}
\caption{The positive braid
$P=\sigma_1\sigma_3\sigma_2\sigma_3$ standardizes
the thick curve $C$ in (a).
The strands in $K''=l_1\cup l_2$ cross $l_3$ once and $l_4$ twice.
The braid diagram $L$ is the union of $K''$ and two parallel copies of $l_3$.
It represents a positive braid $Q$ which standardizes both $C$ and $C'$.}
\label{fig:St_4br}
\end{figure}

\begin{proof}
Let $C'$ be a standard curve which is disjoint from $C$
and encloses the punctures $\{r,r+1,\ldots,r+s\}$.
Because $C$ and $C'$ are disjoint,
$C$ is either inside $C'$ or outside $C'$
as Figure~\ref{fig:St_2cur}.
\begin{figure}
\tabcolsep=10pt
\begin{tabular}{ccc}
\includegraphics[scale=.8]{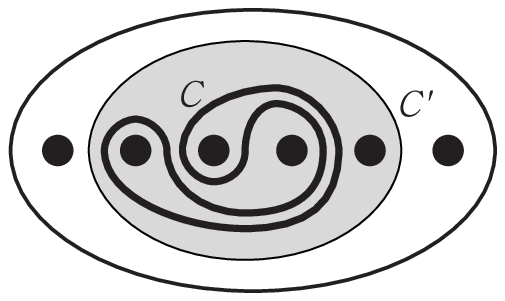}&
\includegraphics[scale=.8]{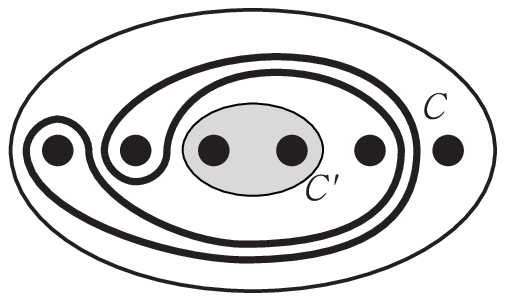}&
\includegraphics[scale=.8]{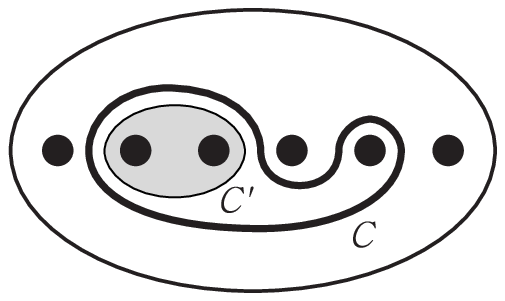}
\end{tabular}
\vskip -2mm
\caption{$C$ is inside $C'$ in the first figure,
and outside $C'$ in the other two figures.}
\label{fig:St_2cur}
\end{figure}

\medskip\noindent
\textsl{Case 1.}\ $C$ is inside $C'$.\\
There exists a positive braid $Q$ written as a positive word on
$\sigma_r,\ldots,\sigma_{r+s-1}$ such that $Q*C$ is standard.
Since $Q\in\St(C)$ and $P$ is the $\le_R$-minimal element of $\St(C)$,
we have $P\le_R Q$, hence $Q=RP$ for some positive braid $R$.
In particular, $P$ is written as a positive word on
$\sigma_r,\ldots,\sigma_{r+s-1}$, and hence $P*C'=C'$ is standard.

\medskip\noindent
\textsl{Case 2.}\ $C$ is outside $C'$.\\
For a braid diagram $K$, let $c(K)$ denote the number of crossings in $K$.
Note that if all the crossings in $K$ are positive, then $K$ represents
a positive braid $Q$ with $|Q|=c(K)$, where $|Q|$ denotes the
word length of $Q$ with respect to $\sigma_i$'s.

\smallskip\noindent
\emph{Claim.\ \
Let $C$ and $C'$ be essential simple closed curves in $D_n$
such that $C'$ is standard and $C$ is outside $C'$.
Let $P$ be an element (not necessarily the $\le_R$-minimal element) of\/ $\St(C)$.
Then there is a positive braid $Q$ such that $|Q|\le |P|$ and
both $Q*C$ and $Q*C'$ are standard.
}

\begin{proof}[Proof of Claim]
See Figure~\ref{fig:St_4br}
which illustrates this proof with a simple example.
Let $K=l_1\cup\cdots\cup l_n$ be a braid diagram
of $P$ in $[0,1]\times\mathbb R$
such that the number of crossings in $K$ is exactly $|P|$.
Here we assume that the right end of $l_i$ is $(1,i)$ for $i=1,\ldots,n$.
Let $\{r,r+1,\ldots,r+s\}$ be the set of punctures inside $C'$.
Let $K'=l_r\cup l_{r+1}\cup\cdots\cup l_{r+s}$
and $K''=K\setminus K'$.
For $i=r,\ldots,r+s$,
let $e_i$ be the number of crossings between $l_i$ and $K''$.
Let $e_{i_0}$ be the minimum of $\{e_r,e_{r+1},\ldots,e_{r+s}\}$.
Then
$$
|P|=c(K)=c(K')+c(K'')+(e_r+\cdots+e_{r+s})
\ge c(K'') + (s+1) e_{i_0}.
$$
Let $L$ be the braid diagram which is the union of $K''$ and
$(s+1)$ parallel copies of $l_{i_0}$,
and let $Q$ be the positive braid represented by $L$.
Since all the crossings in $L$ are positive,
$$
|Q|=c(L) = c(K'')+(s+1)e_{i_0}\le|P|.
$$
By the construction of $Q$, both the curves $Q*C$ and $Q*C'$ are standard.
\end{proof}

By the above claim, there exists a positive braid $Q$ such that
$|Q|\le|P|$ and both $Q*C$ and $Q*C'$ are standard.
Because $P$ is the $\le_R$-minimal element of $\St(C)$ and $Q*C$ is standard,
we have $P\le_R Q$. Since $|Q|\le |P|$, we obtain $P=Q$,
hence $P*C'$ is standard.
\end{proof}

Proposition~\ref{thm:FixStCur} says that
if we standardize the components of a curve system
$\C=C_1\cup\cdots\cup C_k$ one after another
by the $\le_R$-minimal element of the standardizers as follows,
then the product of the $\le_R$-minimal elements used in this process
is exactly the $\le_R$-minimal element of $\St(\C)$.

\begin{enumerate}
\item[(i)]
Standardize the first component $C_1$ of $\C$
using the $\le_R$-minimal element $P_1$ of $\St(C_1)$.
Then $P_1*\C=P_1*C_1\cup \cdots \cup P_1*C_k$ and $P_1*C_1$ is standard.

\item[(ii)]
Standardize the second component $P_1*C_2$ of $P_1*\C$
by the $\le_R$-minimal element $P_2$ of $\St(P_1*C_2)$.
Then the first two components $(P_2P_1)*(C_1\cup C_2)$ of $(P_2P_1)*\C$
are standard.

\item[(iii)]
Continue the above process. Then $(P_k\cdots P_1)*\C$ is standard.
Corollary~\ref{cor:StCur} shows that
in fact $P_k\cdots P_1$ is the $\le_R$-minimal element
of $\St(\C)$.
\end{enumerate}

\begin{corollary}\label{cor:StCur}
Let $\C, \C_1,\ldots, \C_k$ be essential curve systems in $D_n$
such that $\C=\C_1\cup\cdots\cup \C_k$.
Let $P$ be the $\le_R$-minimal element of\/ $\St(\C)$.
\begin{enumerate}
\item[(i)] If $P_i$ is the $\le_R$-minimal element of\/
$\St((P_{i-1}\cdots P_1)*\C_i)$, then $P=P_kP_{k-1}\cdots P_1$.
\item[(ii)] For any standard curve $C'$ disjoint from $\C$,
the curve $P*C'$ is standard.
\end{enumerate}
\end{corollary}

\begin{proof}
We prove the corollary only for the case when
each curve system $\C_i$ has only one component.
The general case can be proved easily from this.
Suppose that each curve system $\C_i$ has only one component.

\medskip\noindent
\emph{Claim.}\ \
The following hold for each $i=0,1,\ldots,k$.
\begin{itemize}
\item[(a)] $P_iP_{i-1}\cdots P_1\le_R P$.
\item[(b)] The curve $(P_iP_{i-1}\cdots P_1)*\C_j$ is standard for $j=1,\ldots,i$.
\item[(c)] For any standard curve $C'$ disjoint from $\C$,
the curve $(P_iP_{i-1}\cdots P_1)*C'$ is standard.
\end{itemize}

\begin{proof}[Proof of Claim]
The statement is obvious for $i=0$ since $P_i\cdots P_1$ is the identity.
Using induction on $i$,
assume that the statement is true for some $i$ with $0\le i<k$.
Since $P_i\cdots P_1\le_R P$,
$$
P=Q(P_i\cdots P_1)
$$
for some $Q\in B_n^+$.
Since $Q*((P_i\cdots P_1)*\C_{i+1})=P*\C_{i+1}$ is standard
and $P_{i+1}$ is the $\le_R$-minimal element of
$\St((P_i\cdots P_1)*\C_{i+1})$, we have
$P_{i+1}\le_R Q$, hence
$$
P_{i+1}P_i\cdots P_1\le_R Q(P_i\cdots P_1)=P .
$$
By the induction hypothesis,
$(P_i\cdots P_1)*C'$ and $(P_i\cdots P_1)*\C_j$
are standard curves disjoint from $(P_i\cdots P_1)*\C_{i+1}$
for $j=1,\ldots,i$.
Since $P_{i+1}$ is the $\le_R$-minimal element
of $\St((P_i\cdots P_1)*\C_{i+1})$,
$(P_{i+1}P_i\cdots P_1)*C'$ and
$(P_{i+1}P_i\cdots P_1)*\C_j$ for $j=1,\ldots,i$
are standard by Proposition~\ref{thm:FixStCur}.
By definition of $P_{i+1}$, $(P_{i+1}P_i\cdots P_1)*\C_{i+1}$ is standard.
\end{proof}

By (b) of the above claim, $(P_kP_{k-1}\cdots P_1)*\C$ is standard.
Since $P$ is the $\le_R$ minimal element of $\St(\C)$,
$P\le_R (P_kP_{k-1}\cdots P_1)$.
By (a) of the claim, $(P_kP_{k-1}\cdots P_1)\le_R P$,
hence $P=P_kP_{k-1}\cdots P_1$.
By (c) of the claim, $P*C'$ is standard for any standard curve $C'$
disjoint from $\C$.
\end{proof}

In the rest of this section, we use the following definition.

\begin{definition}
For a composition $\n=(n_1,\ldots,n_k)$ of $n$,
we define the symbol $\delta_\n$ and non-negative integers
$N_0,N_1,\ldots,N_k$ as follows:
\begin{itemize}
\item $\delta_\n=\Delta_1\oplus\cdots\oplus \Delta_k$,
where $\Delta_i$ is the fundamental braid of $B_{n_i}$ for $i=1,\ldots,k$;
\item $N_0=0$ and $N_i=n_1+n_2+\cdots+n_i$ for $i=1,\ldots,k$.
\end{itemize}
\end{definition}

Then, for a composition $\n=(n_1,\ldots,n_k)$ of $n$
and $\sigma_i\in B_k$, the following hold.

\begin{figure}
\includegraphics[scale=.9]{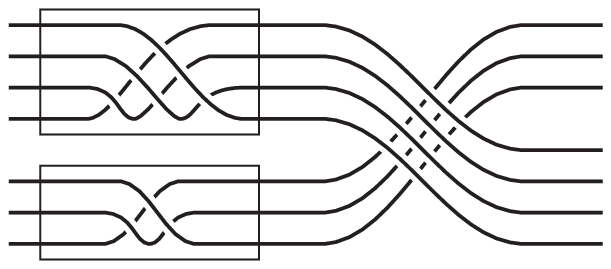}
\quad\raisebox{10mm}{\LARGE $=$}\quad
\includegraphics[scale=.9]{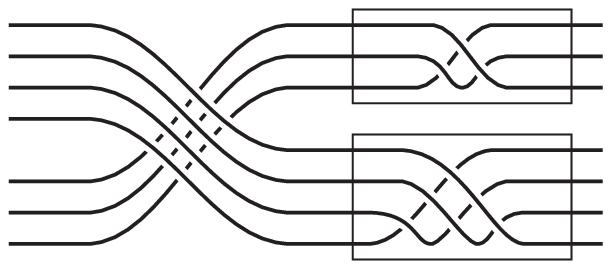}
\caption{The figure shows that
$\delta_{(3,4)}\myangle{\sigma_1}_{(4,3)}=
\myangle{\sigma_1}_{(4,3)} \delta_{(4,3)}$.}
\label{fig:Prod}
\end{figure}

\begin{itemize}
\item
If $A\le_L \delta_\n$, then $A*\C_\n=\C_\n$.
\item
$S(\delta_\n)=F(\delta_\n)
=\{1,\ldots,n-1\}\setminus\{N_1,\ldots,N_{k-1}\}$.
\item
$\sigma_i*\n=\sigma_i^{-1}*\n
=(n_1,\ldots,n_{i-1},n_{i+1},n_i,n_{i+2},\ldots,n_k)$.
\item
$\delta_\n \myangle{\sigma_i}_{\sigma_i*\n}
= \myangle{\sigma_i}_{\sigma_i*\n} \delta_{\sigma_i*\n}$.
See Figure~\ref{fig:Prod}.
\end{itemize}

\begin{lemma}\label{thm:delta_n}
Let $\n=(n_1,\ldots,n_k)$ be a composition of $n$.
\begin{enumerate}
\item[(i)]
Let $A$ be a permutation $n$-braid with induced permutation $\theta$.
Then $\delta_\n A$ is a permutation braid if and only if\/
$\theta^{-1}$ is
order-preserving on the set $\{ N_{i-1}+1,\ldots,N_i\}$
for each $i=1,\ldots,k$,
that is,
$$
\theta^{-1}(N_{i-1}+1)
< \theta^{-1}(N_{i-1}+2)
< \cdots
< \theta^{-1}(N_i).
$$
\item[(ii)]
For a positive $n$-braid $P$,
the starting set $S(\delta_\n P)$ is strictly greater than
the starting set $S(\delta_\n)$
if and only if $\myangle{\sigma_i}_{\sigma_i*\n}\le_L P$
for some $i\in\{1,\ldots, k-1\}$.
\end{enumerate}
\end{lemma}

\begin{proof}
(i)
It is an easy consequence of the fact that
a positive braid $P$ is a permutation braid if and only if
any two of its strands cross at most once~\cite[Lemma9.1.10]{Thu92}
or \cite[Lemma 2.3]{EM94}.
See Figure~\ref{fig:OrePre}.

\begin{figure}
\includegraphics[scale=1]{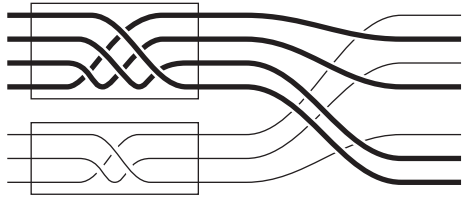}
\caption{The figure shows a permutation braid of the form
$\delta_{(3,4)} A$ for a permutation braid $A$.
If $\theta$ is the induced permutation of $A$, then $\theta^{-1}$ is
order-preserving on each of the sets $\{1,2,3\}$ and $\{4,5,6,7\}$.}
\label{fig:OrePre}
\end{figure}

\smallskip
(ii)
See Figure~\ref{fig:LeftWt}.
\begin{figure}
\includegraphics[scale=1]{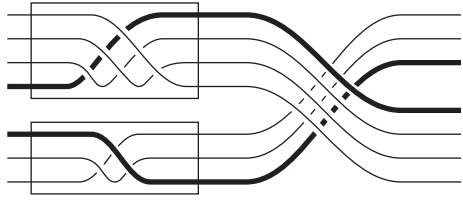}
\caption{The figure shows a permutation braid of the form
$\delta_{(3,4)}A$ for a permutation braid $A$.
If $3\in S(\delta_{(3,4)}A)$, then
two thick strands cross each other and,
hence, $\myangle{\sigma_1}_{(4,3)}\le_L A$.}
\label{fig:LeftWt}
\end{figure}
Suppose $\myangle{\sigma_i}_{\sigma_i*\n}\le_L P$ for some $i\in\{1,\ldots, k-1\}$.
Then $N_i\in S(\delta_\n P)$, hence $S(\delta_\n P)$ is strictly
greater than $S(\delta_\n)$.
Conversely, suppose that $S(\delta_\n P)$ is strictly
greater than $S(\delta_\n)$.
Let $A$ be the permutation $n$-braid such that
$\LM(\delta_\n P)=\delta_\n A$, that is,
$\delta_\n A$ is the first permutation braid in the
left normal form of $\delta_\n P$.
Then $N_i\in S(\delta_\n A)$ for some $i\in\{1,\ldots,k-1\}$.
Let $\omega$ and $\theta$ be the induced permutations of
$\delta_\n$ and $A$ respectively.
Then
$$
\omega^{-1}(N_i)=N_{i-1}+1\quad\mbox{and}\quad
\omega^{-1}(N_i+1)=N_{i+1}.
$$
Since $N_i\in S(\delta_\n A)$, we have
$(\omega\theta)^{-1}(N_i)>(\omega\theta)^{-1}(N_i+1)$ and, hence,
\begin{equation}\label{eqn:St1}
\theta^{-1}(N_{i-1}+1)>\theta^{-1}(N_{i+1}).
\end{equation}
Because $\theta^{-1}$ is order-preserving on each of the sets
$\{N_{i-1}+1,N_{i-1}+2,\ldots,N_i\}$ and
$\{N_i+1,N_i+2,\ldots,N_{i+1}\}$,
we have the following:
\begin{eqnarray}
\theta^{-1}(N_{i-1}+1)<\cdots< \theta^{-1}(N_i-1) < \theta^{-1}(N_i);
\label{eqn:St2}\\
\theta^{-1}(N_i+1)< \theta^{-1}(N_i+2)<\cdots< \theta^{-1}(N_{i+1}).
\label{eqn:St3}
\end{eqnarray}
From (\ref{eqn:St1}), (\ref{eqn:St2}) and (\ref{eqn:St3}),
we obtain $\myangle{\sigma_i}_{\sigma_i*\n}\le_L A\le_L P$.
\end{proof}

The following proposition characterizes the minimal element of
the standardizer $\St(\C)$ of a curve system $\C$.

\begin{proposition}\label{thm:MinElt}
Let $\C$ be an unnested curve system in $D_n$.
Let $P$ be a positive braid such that $P*\C$ is standard and,
hence, $P*\C=\C_\n$ for some composition $\n=(n_1,\ldots,n_k)$ of\/ $n$.
Then the following conditions are equivalent.
\begin{enumerate}
\item[(i)]
$P$ is the $\le_R$-minimal element of the standardizer $\St(\C)$.
\item[(ii)]
$P\wedgeL \delta_\n=1$ and $S(\delta_\n P)=S(\delta_\n)$.
\item[(iii)]
$P^{-1} (\delta_\n P)$ is in \emph{np}-form.
\item[(iv)]
$P^{-1} (\delta_\n^l P)$ is in \emph{np}-form for some $l\ge 1$.
\item[(v)]
$P^{-1} (\delta_\n^l P)$ is in \emph{np}-form for all $l\ge 1$.
\end{enumerate}
\end{proposition}

\begin{proof}
We prove the equivalence by showing that
(i) $\Leftrightarrow$ (ii) $\Rightarrow$
(v) $\Rightarrow$ (iii) $\Rightarrow$ (iv) $\Rightarrow$ (ii).
The implications (v) $\Rightarrow$ (iii) and (iii) $\Rightarrow$ (iv) are obvious.

\smallskip
(i) $\Rightarrow$ (ii)\ \
Let $A=P\wedgeL \delta_\n$ and let $P=AQ$ for some positive braid $Q$.
Since $A\le_L\delta_\n$, $A*\C_\n=\C_\n$, and hence
$$
Q*\C=A^{-1}*(P*\C)=A^{-1}*\C_\n=\C_\n.
$$
Therefore $Q\in\St(\C)$.
By the $\le_R$-minimality of $P$, we have $P=Q$ and, hence,
$P\wedgeL \delta_\n = A=1$.

Assume that $S(\delta_\n P)$ is strictly greater than $S(\delta_\n)$.
Then, by Lemma~\ref{thm:delta_n}~(ii),
$P=\myangle{\sigma_i}_{\sigma_i*\n}Q$ for some $i\in\{ 1,\ldots, k-1\}$
and some positive braid $Q$.
Since
$$
Q*\C
=(\myangle{\sigma_i}_{\sigma_i*\n})^{-1}*(P*\C)
=\myangle{\sigma_i^{-1}}_\n *\C_\n,
$$
$Q*\C$ is standard.
This contradicts the $\le_R$-minimality of $P$.
Consequently, $S(\delta_\n P)= S(\delta_\n)$.

\medskip
(ii) $\Rightarrow$ (i)\ \
Let $Q$ be the $\le_R$-minimal element of $\St(\C)$.
Let $Q*\C=\C_{\n'}$ for some composition $\n'$ of $n$.
Since $P*\C$ is standard, $P=RQ$ for some positive braid $R$.
Since
$$
R*\C_{\n'}=R*(Q*\C)=P*\C=\C_\n,
$$
the positive braid $R$ sends the standard curve system $\C_{\n'}$ to
the standard curve system $\C_\n$.
Therefore, by Lemmas~\ref{thm:decom}~(ii) and~\ref{thm:infofredbr}~(iii),
$R=\myangle{R_0}_{\n'}(R_1\oplus\cdots\oplus R_k)$ for some
positive braids $R_i$ with appropriate braid indices,
and $R_0*\n'=\n$.

If $(R_1\oplus\cdots\oplus R_k)\ne 1$, then
$P\wedgeL \delta_\n\ne 1$.
This contradicts the hypothesis.
Therefore $(R_1\oplus\cdots\oplus R_k)= 1$.

If $R_0\ne 1$, then $R_0=\sigma_i R_0'$ for some $i\in\{ 1,\ldots,k-1\}$
and a positive $k$-braid $R_0'$.
Since $R_0'*\n'=\sigma_i^{-1}*(R_0*\n')=\sigma_i^{-1}*\n=\sigma_i*\n$,
$$
\myangle{R_0}_{\n'}
= \myangle{\sigma_i}_{R_0'*\n'} \myangle{R_0'}_{\n'}
= \myangle{\sigma_i}_{\sigma_i*\n} \myangle{R_0'}_{\n'}.
$$
Since $\myangle{\sigma_i}_{\sigma_i*\n}\le_L \myangle{R_0}_{\n'}\le_LP$,
$S(\delta_\n P)$ is strictly greater than $S(\delta_\n)$
by Lemmas~\ref{thm:decom}~(ii).
This contradicts the hypothesis $S(\delta_\n P)=S(\delta_\n)$.
Therefore $R=1$ and, hence, $P$ is the $\le_R$-minimal
element of $\St(\C)$.

\smallskip
(ii) $\Rightarrow$ (v)\ \
We first claim that $S(\delta_\n^l P)=S(\delta_\n)$ for all $l\ge 1$.
Let $\delta_\n A=\LM(\delta_\n P)$.
Since $S(\delta_\n A)=S(\delta_\n P)$ by Lemma~\ref{lemma:starting_set}~(i)
and $S(\delta_\n P)=S(\delta_\n)$ by the hypothesis,
$$
S(\delta_\n A)=S(\delta_\n P)=S(\delta_\n)=F(\delta_\n).
$$
In particular, $F(\delta_\n)\supset S(\delta_\n A)$,
and hence
$\delta_\n (\delta_\n A)$ is in left normal form
by Lemma~\ref{lemma:starting_set}~(iii).
Since $F(\delta_\n)=S(\delta_\n)$,
$\underbrace{\delta_\n\cdots\delta_\n}_{l-1} (\delta_\n A)$ is the left normal form
of $\delta_\n^l A$ for all $l\ge 1$,
and hence $S(\delta_\n^l P)=S(\delta_\n^l A)=S(\delta_\n)$.

Now we have $S(\delta_\n^l P)=S(\delta_\n)$ for all $l\ge 1$.
By the hypothesis $P\wedgeL\delta_\n=1$,
$$
S(P)\cap S(\delta_\n^l P)=S(P)\cap S(\delta_\n)=\emptyset
\quad\mbox{for all } l\ge 1.
$$
Consequently,
$P\wedgeL \delta_\n^l P=1$ and $P^{-1}(\delta_\n^l P)$ is in
\emph{np}-form for all $l\ge 1$.

\smallskip
(iv) $\Rightarrow$ (ii)\ \
Let $P^{-1}(\delta_\n^l P)$ is in \emph{np}-form for some $l\ge 1$,
that is, $P\wedgeL (\delta_\n^l P)=1$.
Since $P\wedgeL \delta_\n\le_L P\wedgeL (\delta_\n^l P)$,
we have $P\wedgeL \delta_\n=1$.

Assume that $S(\delta_\n P)$ is strictly greater than $S(\delta_\n)$.
By Lemma~\ref{thm:delta_n}~(ii), we have
\begin{equation}\label{eqn:St-Char1}
P=\myangle{\sigma_i}_{\sigma_i*\n}Q
\end{equation}
for some $i\in\{ 1,\ldots, k-1\}$ and some positive braid $Q$.
Since $\delta_\n \myangle{\sigma_i}_{\sigma_i*\n}
= \myangle{\sigma_i}_{\sigma_i*\n} \delta_{\sigma_i*\n}$,
\begin{equation}\label{eqn:St-Char2}
\delta_\n^l P=\delta_\n^l \myangle{\sigma_i}_{\sigma_i*\n}Q
 = \myangle{\sigma_i}_{\sigma_i*\n} \delta_{\sigma_i*\n}^{l}Q.
\end{equation}
By (\ref{eqn:St-Char1}) and (\ref{eqn:St-Char2}),
we obtain $\myangle{\sigma_i}_{\sigma_i*\n} \le_L P\wedgeL (\delta_\n^l P)$,
which contracts the hypothesis that
$P^{-1}(\delta_\n^l P)$ is in \emph{np}-from.
As a result, $S(\delta_\n P)=S(\delta_\n)$.
\end{proof}

Now we are ready to show that standardizing a reduction system $\C$
of a braid by the $\le_R$-minimal element of $\St(\C)$
preserves the membership of the super summit set, ultra summit set and
stable super summit set.
The anonymous referee of this journal pointed out that
our initial proof of the following theorem contains a mistake.
The proof is corrected as suggested by the referee.

\begin{theorem}\label{thm:StMain}
Let $\alpha$ be a reducible $n$-braid with a reduction system $\C$.
Let $P$ be the $\le_R$-minimal element of\/ $\St(\C)$.
Then the following hold.
\begin{enumerate}
\item[(i)] $\inf(\alpha)\le\inf(P\alpha P^{-1})\le\sup(P\alpha P^{-1})\le\sup(\alpha)$.
\item[(ii)] If\/ $\alpha\in[\alpha]^S$, then $P\alpha P^{-1}\in[\alpha]^S$.
\item[(iii)] If\/ $\alpha\in[\alpha]^U$, then $P\alpha P^{-1}\in[\alpha]^U$.
\item[(iv)] If\/ $\alpha\in[\alpha]^\ST$, then $P\alpha P^{-1}\in[\alpha]^\ST$.
\end{enumerate}
\end{theorem}

\begin{proof}
First, suppose that $\C$ is an unnested curve system.
Let $P*\C=\C_\n$ for a composition $\n=(n_1,\ldots,n_k)$ of $n$.
Let $u=\sup(P)$. Define
$\bar P = \Delta^u P^{-1}$ and $Q = \bar P \delta_{\n}^{2} P$.
By Proposition~\ref{thm:MinElt},
$P^{-1}(\delta_\n^2 P)$ is in \emph{np}-form,
hence, by Lemma~\ref{lemma:np-form}~(i),
$$
\bar P=Q\wedgeL\Delta^{\sup(\bar P)}.
$$
Since $(P\alpha P^{-1})*\C_\n=\C_\n$,
$P\alpha P^{-1}=\myangle{\beta_0}_\n(\beta_1\oplus\cdots\oplus \beta_k)$
for some $\beta_i$'s with appropriate braid indices,
and $\beta_0*\n=\n$.
Thus $P\alpha P^{-1}$ commutes with $\delta_\n^{2}$,
and it follows that $\alpha$ commutes with $P^{-1}\delta_\n^{2} P$.
Therefore
$Q\alpha Q^{-1}
= \left( \Delta^u P^{-1} \delta_{\n}^{2} P \right)
\alpha
\left( \Delta^u P^{-1} \delta_{\n}^{2} P \right)^{-1}
= \tau^{-u}(\alpha)$.
That is,
\begin{equation}\label{eqn:st4}
Q^{-1} \tau^{-u}(\alpha) Q = \alpha.
\end{equation}

Consider the following sets:
\begin{eqnarray*}
C(\alpha)&=& \{R\in B_n^+:\inf(\alpha)\le\inf(R^{-1}\alpha R)
\le\sup(R^{-1}\alpha R)\le\sup(\alpha)\};\\
C^S(\alpha)&=& \{R\in B_n^+: R^{-1}\alpha R\in[\alpha]^S\};\\
C^U(\alpha)&=& \{R\in B_n^+: R^{-1}\alpha R\in[\alpha]^U\};\\
C^\ST(\alpha)&=&\{R\in B_n^+: R^{-1}\alpha R\in[\alpha]^\ST\}.
\end{eqnarray*}
By Franco and Gonz\'alez-Meneses~\cite{FG03}, Gebhardt~\cite{Geb05}
and Lee and Lee~\cite{LL06a}, all the sets
$C(\alpha)$, $C^S(\alpha)$, $C^U(\alpha)$ and $C^\ST(\alpha)$ are
closed under $\wedgeL$.

Suppose $\alpha\in[\alpha]^S$.
Since $\tau^m(\alpha)\in [\alpha]^S$ for all $m\in\Z$,
$\Delta^{\sup(\bar P)}\in C^S(\tau^{-u}(\alpha))$.
Since $Q\in C^S(\tau^{-u}(\alpha))$ by (\ref{eqn:st4}),
we have $\bar P=Q\wedgeL\Delta^{\sup(\bar P)} \in C^S(\tau^{-u}(\alpha))$.
That is,
$$
P\alpha P^{-1}
= {\bar P}^{-1} \Delta^u \alpha \Delta^{-u} \bar P
= {\bar P}^{-1} \tau^{-u}(\alpha)\bar P
 \in [ \tau^{-u}(\alpha) ]^S
= [\alpha]^S.
$$
Hence (ii) is proved. The other statements can be proved similarly.
\medskip

Now we consider general case.
For a reduction system $\C$ of $\alpha$,
we decompose $\C$ into $\C_1\cup\cdots \cup \C_l$, where $\C_i$'s are inductively
defined as the outermost component of $\C\setminus(\C_1\cup\cdots\cup\C_{i-1})$.
By the construction,
$\C_1,\ldots, \C_l$ are unnested reduction systems of $\alpha$.
For $i=1,\ldots,l$,
define positive braids $P_i$ and conjugates $\alpha_i$ of $\alpha$
inductively as follows.
Let $P_0 =1$ and $\alpha_0= \alpha$.
\begin{itemize}
\item
$P_i$ is the $\le_R$-minimal element of $\St((P_{i-1}\cdots P_1)*\C_i)$;

\item
$\alpha_i= P_{i}\alpha_{i-1} P_{i}^{-1}
= (P_{i}\cdots P_1)\alpha (P_{i}\cdots P_1)^{-1}$.
\end{itemize}
Note that each $\alpha_i$ is a reducible braid with a reduction system
$(P_{i}\cdots P_1)*\C_{i+1}$
and that $P=P_l\cdots P_1$ by Corollary~\ref{cor:StCur}~(i).

Suppose $\alpha\in[\alpha]^S$.
By the previous discussion on the unnested case,
$P_{i+1}\alpha_i P_{i+1}^{-1}\in[\alpha]^S$
for $i=0,\ldots,l-1$, hence $P\alpha P^{-1}\in[\alpha]^S$.
Therefore (ii) is proved. The other statements can be proved similarly.
\end{proof}

\section{Outermost components of non-periodic reducible braids}
In this section we define the outermost component $\alpha_\ext$ of
a non-periodic reducible braid $\alpha$
using the $\le_R$-minimal element of the standardizer
of the canonical reduction system of $\alpha$,
and study its properties.

Recall the canonical reduction system of mapping classes.
For a reduction system $\C\subset D_n$ of an $n$-braid $\alpha$,
let $D_\C$ be the closure of $D_n\setminus N(\C)$ in
$D_n$, where $N(\C)$ is a regular neighborhood of $\C$. The
restriction of $\alpha$ induces an automorphism on $D_\C$ that is
well-defined up to isotopy. Due to Birman, Lubotzky and
McCarthy~\cite{BLM83} and Ivanov~\cite{Iva92},
for any $n$-braid $\alpha$,
there is a unique
\emph{canonical reduction system} $\R(\alpha)$ with the following
properties.
\begin{enumerate}
\item[(i)]
$\R(\alpha^m)=\R(\alpha)$ for all $m\ne 0$.

\item[(ii)]
$\R(\beta\alpha\beta^{-1})=\beta*\R(\alpha)$ for all $\beta\in B_n$.

\item[(iii)]
The restriction of $\alpha$ to each component of $D_{\R(\alpha)}$ is
either periodic or pseudo-Anosov. A reduction system with this
property is said to be \emph{adequate}.

\item[(iv)]
If $\C$ is an adequate reduction system of $\alpha$,
then $\R(\alpha)\subset\C$.
\end{enumerate}

By the properties of canonical reduction systems,
a braid $\alpha$ is non-periodic reducible if and only if
$\R(\alpha)\ne\emptyset$.
Let $\Rext(\alpha)$ denote the collection
of the outermost components of $\R(\alpha)$.
Then, $\Rext(\alpha)$ is an unnested curve system satisfying
the properties (i) and (ii).
We remark that, while the canonical reduction systems are defined
for the mapping classes of surfaces with genus,
we have to restrict ourselves to the mapping classes of punctured disks
in order to define the outermost component $\Rext(\alpha)$.

\begin{lemma}\label{thm:commute}
Let $\alpha, \beta \in B_n$ with $\R(\alpha)\neq\emptyset$.
If $\alpha\beta=\beta\alpha$,
then $\R(\alpha)$ and $\Rext(\alpha)$ are reduction systems of\/ $\beta$.
\end{lemma}

\begin{proof}
Since $\R(\alpha)=\R(\beta\alpha\beta^{-1})=\beta*\R(\alpha)$ and
$\Rext(\alpha)=\Rext(\beta\alpha\beta^{-1})=\beta*\Rext(\alpha)$,
both $\R(\alpha)$ and $\Rext(\alpha)$ are reduction systems of $\beta$.
\end{proof}

\begin{definition}
Let $\alpha\in B_n$ with $\R(\alpha)\ne\emptyset$.
Let $P$ be the $\le_R$-minimal element of $\St(\Rext(\alpha))$ and
$\beta=P\alpha P^{-1}$.
Since $\Rext(\beta)$ is unnested and standard,
$\Rext(\beta)=\C_\n$ for a composition
$\n=(n_1,\ldots,n_k)$ of $n$, and
$\beta$ has the unique expression
$\beta=\myangle{\beta_0}_\n(\beta_1\oplus\cdots\oplus\beta_k)$
by Lemma~\ref{thm:decom}~(ii).
We define the \emph{outermost component}
$\alpha_\ext$ of $\alpha$ by $\alpha_\ext=\beta_0$.
\end{definition}

In other words, $\alpha_\ext$ is the restriction of $\alpha$ to
the outermost component of $D_n\setminus \Rext(\alpha)$.
This element is a priori defined up to conjugacy, but the use of
the $\le_R$-minimal element $P$ determines the particular element $\beta_0$
to be chosen in the conjugacy class.

\begin{lemma}\label{thm:infscompare}
Let $\alpha$ be an $n$-braid with $\R(\alpha)\ne\emptyset$.
\begin{itemize}
\item[(i)]
If\/ $\beta$ is conjugate to $\alpha$,
then $\beta_\ext$ is conjugate to $\alpha_\ext$.
\item[(ii)]
$(\alpha^m)_\ext=(\alpha_\ext)^m$ for all $m\ne 0$.
\item[(iii)]
$\inf(\alpha)\le\inf(\alpha_\ext)\le\sup(\alpha_\ext)\le\sup(\alpha)$.
\item[(iv)]
$\infs(\alpha)\le\infs(\alpha_\ext)\le\sups(\alpha_\ext)\le\sups(\alpha)$.
\item[(v)]
$\INF(\alpha)\le\INF(\alpha_\ext)\le\SUP(\alpha_\ext)\le\SUP(\alpha)$.
\end{itemize}
\end{lemma}

\begin{proof}
(i) is obvious.
(ii) follows from $\R(\alpha^m)=\R(\alpha)$.
(iii) follows from Lemma~\ref{thm:infofredbr}
and Theorem~\ref{thm:StMain}.

\noindent
(iv)
Choose any $\beta\in [\alpha]^S$.
By (iii), we have
$$
\infs(\alpha)=\inf(\beta)\le\inf(\beta_\ext)
\le\sup(\beta_\ext)\le\sup(\beta)=\sups(\alpha).
$$
Since $\alpha_\ext$ and $\beta_\ext$ are conjugate by (i),
$$
\inf(\beta_\ext)\le\infs(\alpha_\ext)\le\sups(\alpha_\ext)\le\sup(\beta_\ext).
$$
Combining the above two, we obtain
$\infs(\alpha)\le\infs(\alpha_\ext)\le\sups(\alpha_\ext)\le\sups(\alpha)$.

\noindent
(v)
By (ii) and (iii), for all $m\ge 1$,
$$
\inf(\alpha^m)\le\inf((\alpha^m)_\ext)=\inf((\alpha_\ext)^m)
\le\sup((\alpha_\ext)^m)=\sup((\alpha^m)_\ext)\le\sup(\alpha^m).
$$
Therefore,
$$
\frac{\inf(\alpha^m)}{m}\le \frac{\inf((\alpha_\ext)^m)}{m}
\le\frac{\sup((\alpha_\ext)^m)}{m}\le\frac{\sup(\alpha^m)}{m}.
$$
By taking $m\to\infty$, we obtain the desired inequalities for
$\INF(\cdot)$ and $\SUP(\cdot)$.
\end{proof}

\begin{lemma}\label{thm:standardSSS}
Let\/ $\alpha\in B_n$ with $\Rext(\alpha)$ standard.
Then $\Rext(\tau(\alpha))$,
$\Rext(\c_0(\alpha))$ and $\Rext(\d(\alpha))$ are standard.
Moreover,
\begin{enumerate}
\item[(i)]
$\tau(\alpha)_\ext=\tau(\alpha_\ext);$
\item[(ii)]
$\c_0(\alpha)_\ext=\biggl\{\begin{array}{ll}
\alpha_\ext       & \mbox{if\/ $\inf(\alpha_\ext)>\inf(\alpha);$}\\
\c_0(\alpha_\ext) & \mbox{if\/ $\inf(\alpha_\ext)=\inf(\alpha);$}
\end{array}$
\item[(iii)]
$\d(\alpha)_\ext=\biggl\{\begin{array}{ll}
\alpha_\ext       & \mbox{if\/ $\sup(\alpha_\ext)<\sup(\alpha);$}\\
\d(\alpha_\ext) \ \  & \mbox{if\/ $\sup(\alpha_\ext)=\sup(\alpha)$.}
\end{array}$
\end{enumerate}
\end{lemma}

\begin{proof}
$\Rext(\tau(\alpha))=\Rext(\Delta^{-1}\alpha\Delta)=\Delta^{-1}*\Rext(\alpha)$
is obviously standard.
$\Rext(\c_0(\alpha))$ and $\Rext(\d(\alpha))$ are standard
by Corollary~\ref{thm:BNG0}.
Let $\Rext(\alpha)=\C_\n$ for a composition $\n=(n_1,\ldots,n_k)$ of $n$
and $\alpha=\myangle{\alpha_0}_\n(\alpha_1\oplus\cdots\oplus\alpha_k)$.
Let $\Delta_i$ be the fundamental braid of $B_{n_i}$
for $i=1,\ldots,k$ and $\Delta_0$ be the fundamental braid of $B_k$.
Note that $\alpha_0*\n=\n$ and
$$\Delta
=(\Delta_1\oplus\cdots\oplus\Delta_k) \myangle{\Delta_0}_{\Delta_0^{-1}*\n}
=\myangle{\Delta_0}_{\Delta_0^{-1}*\n}(\Delta_k\oplus\cdots\oplus\Delta_1).
$$
Therefore,
\begin{eqnarray*}
\tau(\alpha) &=& \Delta^{-1}\alpha\Delta\\
&=&
\myangle{\Delta_0^{-1}}_{\n}(\Delta_1^{-1}\oplus\cdots\oplus\Delta_k^{-1})\
\myangle{\alpha_0}_\n(\alpha_1\oplus\cdots\oplus\alpha_k)\
\myangle{\Delta_{0}}_{\Delta_0^{-1}*\n}
  (\Delta_k\oplus\cdots\oplus\Delta_1)\\
&=&
\myangle{\Delta_0^{-1}\alpha_0\Delta_0}_{\Delta_0^{-1}*\n}
  (\Delta_k^{-1}\alpha_k\Delta_k\oplus\cdots
  \oplus\Delta_1^{-1}\alpha_1\Delta_1)\\
&=&
\myangle{\tau(\alpha_0)}_{\Delta_0^{-1}*\n}
  (\tau(\alpha_k)\oplus\cdots
  \oplus\tau(\alpha_1)).
\end{eqnarray*}
Since $\R_\ext(\tau(\alpha))=\Delta^{-1}*\C_\n=\C_{\Delta_0^{-1}*\n}$,
$\tau(\alpha)_\ext=\tau(\alpha_0)=\tau(\alpha_\ext)$.

Let $\alpha=\Delta^uA_1\ldots A_l$ be the left normal form of $\alpha$.
Since $\alpha*\C_\n=\C_\n$ is standard,
$A_l*\C_\n$ is standard by Theorem~\ref{thm:normalNstandard}.
By Lemmas~\ref{thm:decom} (ii) and~\ref{thm:infofredbr}~(iii),
$A_l$ is expressed as
$A_l=\myangle{A_{l,0}}_\n(A_{l,1}\oplus\cdots\oplus A_{l,k})$,
where $A_{l,i}$'s are permutation $n_i$-braids.
Let $\theta_1$ and $\theta_2$ be the induced permutations
of $\alpha_0A_{l,0}^{-1}$ and $A_{l,0}^{-1}$ respectively.
Then
\begin{eqnarray*}
\d(\alpha) &=& A_l\alpha A_l^{-1}\\
&=&
\myangle{A_{l,0}}_\n(A_{l,1}\oplus\cdots\oplus A_{l,k})\
\myangle{\alpha_0}_\n(\alpha_1\oplus\cdots\oplus\alpha_k)\
(A_{l,1}^{-1}\oplus\cdots\oplus A_{l,k}^{-1})
\myangle{A_{l,0}^{-1}}_{A_{l,0}*\n}\\
&=&
\myangle{A_{l,0}\alpha_0A_{l,0}^{-1}}_{A_{l,0}*\n}
(A_{l,\theta_1(1)}\alpha_{\theta_2(1)} A_{l,\theta_2(1)}^{-1}
\oplus\cdots\oplus
A_{l,\theta_1(k)}\alpha_{\theta_2(k)} A_{l,\theta_2(k)}^{-1}).
\end{eqnarray*}
Recall Lemma~\ref{thm:infofredbr}~(ii) that $\sup(\alpha_\ext)\le\sup(\alpha)$.
If $\sup(\alpha_\ext)<\sup(\alpha)$,
then $A_{l,0}=1$ and, hence, $\d(\alpha)_\ext=\alpha_\ext$.
If $\sup(\alpha_\ext)=\sup(\alpha)$,
then $A_{l,0}\neq 1$ and, hence,
$\d(\alpha)_\ext=A_{l,0}\alpha_0A_{l,0}^{-1}=\d(\alpha_\ext)$.

For $\c_0(\alpha)$, use the identity $\c_0(\alpha)=\d(\alpha^{-1})^{-1}$.
\end{proof}

Recall Lemma~\ref{thm:infscompare} that
$\inf(\alpha)\le\inf(\alpha_\ext)$ and
$\infs(\alpha)\le\infs(\alpha_\ext)$
for any $\alpha\in B_n$ with $\R(\alpha)\ne\emptyset$.

\begin{lemma}\label{thm:inf_out}
Let $\alpha$ be an $n$-braid with $\R(\alpha)\ne\emptyset$.
Let $\beta$ be an element of\/ $[\alpha]^U$ with $\Rext(\beta)$ standard.
\begin{itemize}
\item[(i)]
Let $\infs(\alpha_\ext)>\infs(\alpha)$.
Then, $\inf(\beta_\ext)>\inf(\beta)$.

\item[(ii)]
Let $\infs(\alpha_\ext)=\infs(\alpha)$.
Then, $\inf(\beta_\ext)=\inf(\beta)$, and
$\c_0^m(\beta_\ext)=\beta_\ext$ for some $m\ge 1$.
\end{itemize}
\end{lemma}

\begin{proof}
By Lemma~\ref{thm:infscompare}~(i), $\beta_\ext$ and $\alpha_\ext$
are conjugate, hence $\inf(\beta_\ext)\le\infs(\alpha_\ext)$.

We first prove the following claim.

\medskip\noindent
\emph{Claim}.\ \
Let $\inf(\beta_\ext)=\inf(\beta)$.
Then, $\c_0^m(\beta_\ext)=\beta_\ext$ for some $m\ge 1$, and
$\infs(\alpha_\ext)=\inf(\beta_\ext)=\inf(\beta)=\infs(\alpha)$.

\begin{proof}[Proof of Claim]
By Lemma~\ref{thm:standardSSS} (ii), the sequence
$\{\inf(\c_0^i(\beta)_\ext)\}_{i=0}^{\infty}$ is non-decreasing.
Since $\beta\in[\alpha]^U$, $\c_0^m(\beta)=\beta$ for some $m\ge 1$.
Therefore,
$$\inf(\c_0^i(\beta)_\ext)=\inf(\beta_\ext)
\qquad\mbox{for all $i\ge 0$}.
$$
Since $\c_0^i(\beta)\in[\alpha]^U$ for all $i\ge 0$,
we have $\inf(\c_0^i(\beta))=\infs(\alpha)=\inf(\beta)$ for all $i\ge 0$.
Hence
$$
\inf(\c_0^i(\beta)_\ext)=\inf(\beta_\ext)
=\inf(\beta)=\inf(\c_0^i(\beta))
\qquad\mbox{for all $i\ge 0$}.
$$
By Lemma~\ref{thm:standardSSS} (ii),
$$\c_0^i(\beta)_\ext=\c_0^i(\beta_\ext)
\qquad\mbox{for all $i\ge 0$}.
$$
Since $\c_0^m(\beta)=\beta$,
we obtain $\c_0^m(\beta_\ext)=\c_0^m(\beta)_\ext=\beta_\ext$.

By Theorem~\ref{thm:ConjAlgorithm}~(i),
$\inf(\beta_\ext) = \infs(\beta_\ext) = \infs(\alpha_\ext)$.
Therefore, $\infs(\alpha_\ext)=\inf(\beta_\ext)=\inf(\beta)=\infs(\alpha)$.
\end{proof}

(i)
Assume $\inf(\beta_\ext)=\inf(\beta)$.
Then $\infs(\alpha_\ext)=\inf(\beta_\ext)=\inf(\beta)=\infs(\alpha)$
by the above claim.
This contradicts the hypothesis that $\infs(\alpha_\ext)>\infs(\alpha)$,
hence $\inf(\beta_\ext)>\inf(\beta)$.

\medskip
(ii)
Since $\inf(\beta)\le\inf(\beta_\ext)\le\infs(\alpha_\ext)$,
$$
\inf(\beta)\le\inf(\beta_\ext)
\le\infs(\alpha_\ext)=\infs(\alpha)=\inf(\beta).
$$
Therefore $\inf(\beta_\ext)=\inf(\beta)$.
By the claim, $\c_0^m(\beta_\ext)=\beta_\ext$ for some $m\ge 1$.
\end{proof}

The following proposition show that the property
$\infs(\alpha_\ext)>\infs(\alpha)$ is preserved by
taking powers.

\begin{proposition}\label{thm:infs-power}
If\/ $\infs(\alpha_\ext)>\infs(\alpha)$, then
$\infs((\alpha^m)_\ext)>\infs(\alpha^m)$
for all $m\ge 1$.
\end{proposition}

\begin{proof}
By Theorem~6.1 in~\cite{Lee07}, for any $\beta\in B_n$ and any $m\ge 1$,
$$\infs(\beta)\le\frac{\infs(\beta^m)}m<\infs(\beta)+1.
$$
By Lemma~\ref{thm:infscompare}~(ii),
$(\alpha^m)_\ext=(\alpha_\ext)^m$ for all $m\ne 0$.
Since $\infs(\alpha_\ext)>\infs(\alpha)$,
$$
\frac{\infs(\alpha^m)}m<\infs(\alpha)+1
\le\infs(\alpha_\ext)\le\frac{\infs((\alpha_\ext)^m)}m
=\frac{\infs((\alpha^m)_\ext)}m
$$
for all $m\ge 1$.
Therefore $\infs((\alpha^m)_\ext)>\infs(\alpha^m)$
for all $m\ge 1$.
\end{proof}

\section{Split braids}
An $n$-braid $\alpha$ is called a \emph{split braid}
if it is conjugate to an element in the subgroup of $B_n$ generated
by $\sigma_1,\ldots,\sigma_{i-1},\sigma_{i+1},\ldots,\sigma_{n-1}$
for some $1\le i\le n-1$~\cite{Hum91}.
In our terminology, $\alpha\in B_n$ is a split braid
if it is conjugate to a braid $\beta$ of the form
$\beta=\myangle{1}_\n(\beta_1\oplus\beta_2)$,
where $\n=(i,n-i)$ for some $1\le i\le n-1$,
and $\beta_1\in B_i$ and $\beta_2\in B_{n-i}$.

The following lemma is easy to show, but we include a proof for completeness.

\begin{lemma}\label{lem:split-equiv}
Let $\alpha$ be an $n$-braid.
\begin{itemize}
\item[(i)]
$\alpha$ is a split braid if and only if either $\alpha$ is the identity
or $\alpha$ is non-periodic and reducible with $\alpha_\ext =1$.

\item[(ii)]
Let\/ $\alpha=\myangle{1}_\n(\alpha_1\oplus\cdots\oplus\alpha_k)$
for a composition $\n=(n_1,\ldots,n_k)$ of\/ $n$
such that\/ $\Rext(\alpha)\neq\emptyset$.
Then $\Rext(\alpha)=\C_\n$ if and only if
$\alpha_i$ is non-split for each $1\le i\le k$.
\end{itemize}
\end{lemma}

\begin{proof}
For unnested curve systems $\C$ and $\C'$ in $D_n$,
let us write ``$\C'\preceq \C$'' if each component of $\C'$ is enclosed by
(possibly parallel to) a component of $\C$,
and write ``$\C'\precneqq \C$'' if $\C'\preceq \C$ and $\C'\ne \C$.
Then $\preceq$ is a partial order over the set of unnested curve systems in $D_n$.
For compositions $\n=(n_1,\ldots,n_k)$ and $\n'$ of $n$,
$\C_{\n'}\preceq\C_\n$ if and only if $\n'$ is a refinement of $\n$,
that is, for each $1\le i\le k$, there exists a composition
$(n'_{i,1}, \ldots, n'_{i, r_i})$ of $n_i$ such that
$\n'=(n_{1,1}',\ldots,n_{1, r_1}',\ldots,
n_{k, 1}',\ldots,n_{k, r_k}')$.

\medskip\noindent
\emph{Claim}.\ \
Let $\beta*\C_\n=\C_\n$ for a composition $\n$ of $n$,
then $\beta$ is written as
$\beta=\myangle{\beta_0}_\n(\beta_1\oplus\cdots\oplus \beta_k)$.
If $\beta_0$ is periodic or pseudo-Anosov, then
$\Rext(\beta)\preceq \C_\n$, and there exists $P\in B_n^+$ such that
both $P*\Rext(\beta)$ and $P*\C_\n$ are standard.

\begin{proof}[Proof of Claim]
Because $\beta_0$ is periodic or pseudo-Anosov,
we can make an adequate reduction system of $\beta$ from $\C_\n$
by adding some curves each of which is enclosed by a component of $\C_\n$.
Because any adequate reduction system of $\beta$
contains $\Rext(\beta)$ as a subset,
we have $\Rext(\beta)\preceq\C_\n$.
Let $P$ be the $\le_R$-minimal element of $\St(\Rext(\beta))$.
Then $P*\Rext(\beta)$ is standard by the construction.
Apply Corollary~\ref{cor:StCur} to $\C_\n\backslash\Rext(\beta)$.
Then $P* (\C_\n\backslash\Rext(\beta))$ and hence $P*\C_\n$ are standard.
\end{proof}

(i)\ \
It is obvious that if $\alpha$ is the identity
or $\alpha$ is non-periodic and reducible with $\alpha_\ext=1$
then $\alpha$ is a split braid.
Conversely, suppose that $\alpha$ is a split braid.
Taking a conjugate of $\alpha$ if necessary,
we may assume that
$$
\alpha=\myangle{1}_\n(\alpha_1\oplus\alpha_2),
$$
where $\n=(\ell,n-\ell)$ for some $1\le \ell\le n-1$,
and $\alpha_1\in B_\ell$ and $\alpha_2\in B_{n-\ell}$.

First, assume that $\Rext(\alpha)=\emptyset$,
that is, $\alpha$ is periodic or pseudo-Anosov.
Since split braids are a special type of reducible braids
and since pseudo-Anosov braids cannot be reducible~\cite{FLP79},
$\alpha$ is periodic. Therefore $\alpha^p=\Delta^{2m}$ for some
$p\ne 0$ and $m\in\Z$, hence
$$
\myangle{1}_\n(\alpha_1^p \oplus\alpha_2^p)
=\alpha^p
=\Delta^{2m}
=\myangle{\Delta_{0}^{2m}}_\n
(\Delta_{1}^{2m} \oplus\Delta_{2}^{2m}),
$$
where $\Delta_0$, $\Delta_1$ and $\Delta_2$ are the fundamental braids of
$B_2$, $B_{\ell}$ and $B_{n-\ell}$, respectively.
By Lemma~\ref{thm:decom}~(i), we have $\Delta_{0}^{2m}=1$,
hence $m=0$, and it follows that $\alpha^p=1$.
Because braid groups are torsion-free~\cite{Deh98},
$\alpha$ is the identity.

\def\Out{\operatorname{Out}}

Now, assume that $\Rext(\alpha)\ne\emptyset$, that is,
$\alpha$ is non-periodic and reducible.
For a curve system $\C$ in a punctured disk $D_m$,
let $\Out(D_m\setminus\C)$ denote the outermost component of $D_m\setminus\C$.
By the above claim,
we have $\Rext(\alpha)\preceq \C_\n$ and hence $\C_\n\subset\Out(D_n\setminus\Rext(\alpha))$,
and we may assume that $\Rext(\alpha)$ is standard.
Let $\alpha_\ext$ be a $k$-braid.
Because $\Rext(\alpha)$ is standard,
$\Out(D_n\setminus\Rext(\alpha))$ is canonically diffeomorphic to $D_k$.
Let $\C'$ be the image of $\C_\n$ under this diffeomorphism.
Then $\C'$ is a reduction system of $\alpha_\ext$ such that
the restriction of $\alpha_\ext$ to $\Out(D_k\setminus \C')$ is
the same as the restriction of $\alpha$ to $\Out(D_n\setminus\C_\n)$
which is the identity.
This means that  $\alpha_\ext$ is a split braid.
Because $\alpha_\ext$ is either periodic or pseudo-Anosov,
the discussion in the above paragraph shows that $\alpha_\ext$ is the identity.

\medskip
(ii)\ \
Assume that $\alpha_\ell$ is a split braid for some $1\le \ell\le k$,
hence $\alpha_\ell$ is conjugate to
$\myangle{1}_{\n_\ell}(\alpha_\ell'\oplus\alpha_\ell'')$,
where $\n_\ell =(n'_\ell, n''_\ell)$ is a composition of $n_\ell$,
and $\alpha_\ell'\in B_{n_\ell'}$ and $\alpha_\ell''\in B_{n_\ell''}$.
By taking a conjugate of $\alpha$ if necessary, we may assume that
$$
\alpha=\myangle{1}_{\n'}(\alpha_1\oplus\cdots\alpha_{\ell-1}\oplus
\alpha_\ell'\oplus\alpha_\ell''\oplus\alpha_{\ell+1}\oplus\cdots\oplus\alpha_k),
$$
where $\n'=(n_1,\ldots,n_{\ell-1},n_\ell',n_\ell'',n_{\ell+1},\ldots,n_k)$.
Note that $\C_{\n'} \precneqq \C_{\n}$.
By the claim, $\Rext(\alpha)\preceq\C_{\n'}\precneqq \C_{\n}$,
hence $\C_\n\ne \Rext(\alpha)$.

Conversely, assume that $\Rext(\alpha)\ne\C_\n$.
By the claim, we may assume that
$\Rext(\alpha)\precneqq \C_\n$ and $\Rext(\alpha)$ is standard.
Let $\Rext(\alpha)=\C_{\n'}$ for a composition $\n'$ of $n$.
Then $\n'$ is a refinement of $\n$, hence, for each $i$, there exists a composition
$(n'_{i,1}, \ldots, n'_{i, r_i})$ of $n_i$ such that
$\n'=(n_{1,1}',\ldots,n_{1, r_1}',\ldots,n_{k, 1}',\ldots,n_{k, r_k}')$.
Because $\alpha_\ext$ is the identity by (i), $\alpha$ is written as
$$
\alpha=\myangle{1}_{\n'}(
\underbrace{\alpha_{1,1}\oplus\cdots\oplus\alpha_{1,r_1}}_{r_1}\oplus\cdots\oplus
\underbrace{\alpha_{k,1}\oplus\cdots\oplus\alpha_{k,r_k}}_{r_k}).
$$
Since $\C_{\n'}=\Rext(\alpha)\ne\C_\n$ by the assumption,
we have $r_{\ell} \ge 2$ for some $1\le \ell\le k$.
Comparing the above expression with
$\alpha=\myangle{1}_\n(\alpha_1\oplus\cdots\oplus\alpha_k)$,
we have $\alpha_{\ell}=
\myangle{1}_{\n'_\ell} (\alpha_{\ell, 1}\oplus\cdots\oplus\alpha_{\ell, r_{\ell}})$
where $\n'_\ell =(n'_{\ell, 1}, \ldots, n'_{\ell, r_\ell})$.
Since $r_{\ell}\ge 2$, $\alpha_{\ell}$ is a split braid.
\end{proof}

For $\alpha\in B_n$, let $|\alpha|$
denote the minimal word length of $\alpha$ with respect
to $\{\sigma_1^{\pm 1},\ldots,\sigma_{n-1}^{\pm 1}\}$.
Then $|\alpha|$ is the minimum number of crossings
in the braid diagram of $\alpha$.

\begin{proposition}\label{thm:splitMinimal}
If\/ $\alpha\ne 1$ is a split braid and\/
$|\alpha|$ is minimal in the conjugacy class of $\alpha$,
then $\Rext(\alpha)$ is standard.
\end{proposition}

\begin{proof}
There exists a braid $\beta$ in the conjugacy class of $\alpha$
such that $\Rext(\beta)$ is standard.
Therefore $\Rext(\beta)=\C_\n$
for some composition $\n=(n_1,\ldots,n_k)$ of $n$.
Then by Lemma~\ref{lem:split-equiv}
$$\beta= \myangle{1}_{\n}(\beta_1\oplus\cdots\oplus\beta_k)$$
for some non-split $n_i$-braids $\beta_i$.
We may choose $\beta$ so that $|\beta_i|$ is minimal in the conjugacy class
of $\beta_i$ for each $i\in\{1,\ldots,k\}$.

Since $\alpha$ and $\beta$ are conjugate,
$\alpha=\gamma\beta\gamma^{-1}$ for some $\gamma\in B_n$.
Let $\theta$ be the induced permutation of $\gamma$.
For $i=1,\ldots,k$, let $S_i=\{j: n_1+\cdots+n_{i-1}<j\le n_1+\cdots+n_i\}$
and $T_i=\{\theta(j): j\in S_i\}$.
Let $\gamma_i$ be the result of forgetting the $j$-th strand
from $\gamma$ for all $j\not\in S_i$.
(The strands of a braid are numbered from bottom to top
at its right end.)
See Figure~\ref{fig:split}.
Let $\alpha_i$ be the result of forgetting the $j$-th strand
from $\alpha$ for all $j\not\in T_i$.
Then $\alpha_i=\gamma_i\beta_i\gamma_i^{-1}$ for all $i=1,\ldots,k$.

\begin{figure}
\includegraphics[scale=1]{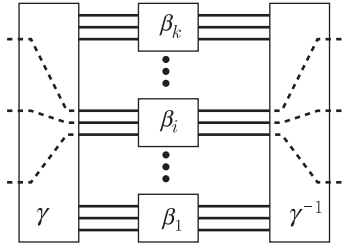}
\caption{The dotted strands indicate $\gamma_i$ in the proof
of Proposition~\ref{thm:splitMinimal}.
}\label{fig:split}
\end{figure}

Let $K$ be a braid diagram of $\alpha$ such that the number of crossings
in $K$ is exactly $|\alpha|$.
For $i=1,\ldots,k$, let $K_i$ be the result of deleting the $j$-th strand
from $K$ for all $j\not\in T_i$.
Then $K_i$ is a braid diagram of $\alpha_i$ for all $i$.
Let $c(K)$ and $c(K_i)$ denote the numbers of crossings in $K$ and $K_i$,
respectively.
Then $|\alpha|=c(K)$,
$|\alpha_i|\le c(K_i)$ for each $i$
and $\sum_{i=1}^k c(K_i)\le c(K)$.

Since $|\alpha|$ is minimal in the conjugacy class, $|\alpha|\le|\beta|$.
Since $|\beta_i|$ is minimal in the conjugacy class,
$|\beta_i|\le |\alpha_i|$ for all $i=1,\ldots,k$.
Hence
$$
c(K)=|\alpha| \le |\beta| = \sum_{i=1}^k|\beta_i|
\le\sum_{i=1}^k|\alpha_i|
\le\sum_{i=1}^k c(K_i)
\le c(K).
$$
Therefore $c(K)=\sum_{i=1}^{k} c(K_i)$ and it follows that
there is no crossing between
the strands in $K_i$ and those in $K_j$ whenever $i\ne j$.

Now we claim that each $T_l$ is a set of consecutive integers.
On the contrary, assume that there exists $j\in T_m$ for some $m\ne l$
such that $i_1<j<i_2$ for some $i_1,i_2\in T_l$.
Let $K_{l,1}$ be the result of deleting all $i$-th strands from $K_l$
with $i>j$ and let $K_{l,2}=K_l\setminus K_{l,1}$.
See Figure~\ref{fig:sp_diagram}.
Because there is no crossing between the strands in $K_l$ and those in $K_m$,
there is no crossing between $K_{l,1}$ and $K_{l,2}$.
Therefore $K_l$ is splitted into $K_{l,1}$ and $K_{l,2}$.
This contradicts that $\alpha_l$ is non-split.
Hence, each $T_l$ is a set of consecutive integers.

\begin{figure}
\includegraphics{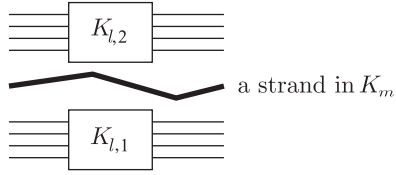}
\caption{Since there is no crossing between the strands in $K_l$ and those in $K_m$,
if a strand of $K_m$ goes through $K_l$, then
$K_l$ is splitted into two parts $K_{l,1}$ and $K_{l,2}$.}
\label{fig:sp_diagram}
\end{figure}

Let $T_{i_1}, T_{i_2}, \ldots, T_{i_k}$ be the rearrangement of
$T_j$'s such that the elements of the sets are increasing,
and let $\n'=( n_{i_1},\ldots, n_{i_k} )$.
Then $\alpha= \myangle{1}_{\n'}(\alpha_{i_1}\oplus\cdots\oplus\alpha_{i_k})$
and $\Rext(\alpha) =\C_{\n'}$.
Therefore $\Rext(\alpha)$ is standard.
\end{proof}

\begin{corollary}\label{thm:PosSplt}
If\/ $P\ne 1$ is a positive split braid, then $\Rext(P)$ is standard.
\end{corollary}

\begin{proof}
If $P$ is a positive braid, then
$|P|$ is minimal in the conjugacy class of $P$.
\end{proof}

\section{Ultra summit sets of reducible braids}
In this section, we establish Theorem~\ref{thm:main},
the main result of this paper.
Roughly speaking, it says that if the outermost component $\alpha_\ext$
is simpler than the whole braid $\alpha$ from a Garside-theoretic point of view,
then it is easy to find a reduction system of $\alpha$.

\begin{definition}
Let $\alpha\in B_n$, $\beta\in [\alpha]^U$ and
$m=\min\{l\ge 1:\c_0^l(\beta)=\beta\}$.
For $i=0,\ldots,m-1$, let $A_i$ be the $\le_R$-minimal element of
$\{P\in B_n^+: \inf(P\c_0^i(\beta))>\inf(\c_0^i(\beta))\}$.
The product $A_{m-1}A_{m-2}\cdots A_0$
is called the \emph{cycling commutator} of $\beta$ and denoted $T_\beta$.
\end{definition}

By definition, the cycling commutator $T_\beta$ is a positive braid.
By Lemma~\ref{thm:duality}~(i),
\begin{eqnarray*}
T_\beta\beta T_\beta^{-1}
&=& A_{m-1}\cdots A_2A_1A_0\ \beta\ A_0^{-1}A_1^{-1}A_2^{-1}\cdots A_{m-1}^{-1}\\
&=& A_{m-1}\cdots A_2A_1\  \c_0(\beta)\ A_1^{-1}A_2^{-1}\cdots A_{m-1}^{-1}\\
&=& A_{m-1}\cdots A_2\ \c_0^2(\beta)\ A_2^{-1}\cdots A_{m-1}^{-1}\\
&=& \cdots =\c_0^m(\beta)=\beta.
\end{eqnarray*}

\begin{lemma}\label{thm:cycl_comm}
Let $\alpha\in B_n$ and $\beta\in[\alpha]^U$.
Then the cycling commutator $T_\beta$ is
a non-identity positive braid with $T_\beta\beta=\beta T_\beta$.
\end{lemma}

The following proposition is the key to Theorem~\ref{thm:main}.
We prove it in~\S\ref{sec:proof}.

\begin{proposition}\label{thm:main-prop}
Let $\alpha$ be a non-periodic reducible $n$-braid
with $\infs(\alpha_\ext)>\infs(\alpha)$.
For any element $\beta$ of\/ $[\alpha]^U$, the cycling commutator $T_\beta$
is a split braid.
\end{proposition}

Recall from Lemma~\ref{thm:infscompare} that
$\infs(\alpha)\le\infs(\alpha_\ext)\le\sups(\alpha_\ext)\le\sups(\alpha)$
and
$\INF(\alpha)\le\INF(\alpha_\ext)\le\SUP(\alpha_\ext)\le\SUP(\alpha)$
for any non-periodic reducible braid $\alpha$.

\begin{theorem}\label{thm:main}
Let $\alpha$ be a non-periodic reducible $n$-braid.
\begin{enumerate}
\item[(i)]
If\/ $\infs(\alpha_\ext)>\infs(\alpha)$,
then each element of\/ $[\alpha]^U$ has a standard reduction system.

\item[(ii)]
If\/ $\sups(\alpha_\ext)<\sups(\alpha)$,
then each element of\/ $[\alpha]_\d^U$ has a standard reduction system.

\item[(iii)]
If\/ $\alpha$ is a split braid,
then each element of\/ $[\alpha]^U \cup [\alpha]_\d^U$ has a standard reduction system.

\item[(iv)]
If\/ $\alpha_\ext$ is periodic,
then there exists $1\le q< n$ such that each element
of\/ $[\alpha^q]^U \cup [\alpha^q]_\d^U$ has a standard reduction system.

\item[(v)]
If\/ $\INF(\alpha_\ext)>\INF(\alpha)$,
then there exists $1\le q< n(n-1)/2$
such that each element
of\/ $[\alpha^q]^U$ has a standard reduction system.

\item[(vi)]
If\/ $\SUP(\alpha_\ext)<\SUP(\alpha)$,
then there exists $1\le q< n(n-1)/2$
such that each element of\/ $[\alpha^q]^U_\d$
has a standard reduction system.
\end{enumerate}
\end{theorem}

\begin{proof}
(i)\
Let $\beta$ be an element of $[\alpha]^U$.
By Proposition~\ref{thm:main-prop}, the cycling commutator
$T_\beta$ is a non-identity positive split braid.
By Corollary~\ref{thm:PosSplt}, $\Rext(T_\beta)$ is standard.
Since $\beta$ commutes with $T_\beta$ by Lemma~\ref{thm:cycl_comm},
$\Rext(T_\beta)$ is
a standard reduction system of $\beta$ by Lemma~\ref{thm:commute}.

\medskip
(ii)\
Because $\infs((\alpha^{-1})_\ext)=\infs((\alpha_\ext)^{-1})=-\sups(\alpha_\ext)$
and $\infs(\alpha^{-1})=-\sups(\alpha)$,
we have $\infs((\alpha^{-1})_\ext)>\infs(\alpha^{-1})$.
By (i), each element of $[\alpha^{-1}]^U$ has a standard reduction system.
Because $[\alpha]^U_\d=\{\beta^{-1}:\beta\in [\alpha^{-1}]^U\}$, we are done.

\medskip
(iii)\
Let $\beta\in[\alpha]^U$.
If $\infs(\alpha)<\infs(\alpha_\ext)$,
then $\beta$ has a standard reduction system
by~(i).
If $\infs(\alpha)=\infs(\alpha_\ext)$, then
$\inf(\beta)=\infs(\alpha)=\infs(\alpha_\ext)=0$ and, hence,
$\beta$ is positive.
Since $\beta$ is split, $\Rext(\beta)$ is standard by Corollary~\ref{thm:PosSplt}.

Since $\alpha$ is a split braid, so is $\alpha^{-1}$.
Thus, every element of $[\alpha^{-1}]^U$ and, hence, $[\alpha]_\d^U$
has a standard reduction system.

\medskip
(iv)\
Let $k$ be the braid index of $\alpha_\ext$.
Because $\alpha_\ext$ is periodic, there exist $1\le q\le k$
and $l\in\Z$ such that
$$
(\alpha_\ext)^q=\Delta_0^{2l},
$$
where $\Delta_0$ is the fundamental braid of $B_k$.
Then $\Delta^{-2l}\alpha^q \neq 1$ is a split braid.
By (iii), every element of
$[\Delta^{-2l}\alpha^q]^U\cup [\Delta^{-2l}\alpha^q]_\d^U$
has a standard reduction system.
Since
$$
[\alpha^q]^U=\{\Delta^{2l}\beta: \beta\in [\Delta^{-2l}\alpha^q]^U\}
\quad\mbox{and}\quad
[\alpha^q]_\d^U=\{\Delta^{2l}\beta: \beta\in [\Delta^{-2l}\alpha^q]_\d^U\},
$$
each element of $[\alpha^q]^U\cup [\alpha^q]_\d^U$ has a standard reduction system.

\medskip
(v)\
Recall from Theorem~\ref{thm:INF} that, for any $\gamma\in B_n$,
\begin{itemize}
\item
$\INF(\gamma)$ is rational with denominator less than
or equal to $|\Delta|=n(n-1)/2$;
\item
$\infs(\gamma)\le\INF(\gamma)<\infs(\gamma)+1$;
\item
$\INF(\gamma^m)=m\INF(\gamma)$ for all integers $m\ge 1$.
\end{itemize}

Let $k$ be the braid index of $\alpha_\ext$.
Then $\INF(\alpha_\ext)=p/q$ for some integers $p$, $q$
with $1\le q\le k(k-1)/2$.
Since $\INF((\alpha_\ext)^q)=q\INF(\alpha_\ext)$ is an integer,
we have $\infs((\alpha_\ext)^q)=q\INF(\alpha_\ext)$. Therefore,
$$
\infs((\alpha^q)_\ext)
=\infs((\alpha_\ext)^q)
=q\INF(\alpha_\ext)
> q\INF(\alpha)
=\INF(\alpha^q)
\ge \infs(\alpha^q).
$$
By (i), every element of $[\alpha^q]^U$ has a standard reduction system.

\medskip
(vi)\
It can be proved in a way similar to (v).
\end{proof}

Now, let us consider the following algorithm.
Let $\alpha$ be a given non-periodic $n$-braid.
\begin{description}
\item[Step 1]
Applying cyclings and decyclings to $\alpha$,
obtain an element $\beta$ of the set $[\alpha]^U\cap[\alpha]^U_\d$
together with an element $\gamma$ such that $\alpha=\gamma\beta\gamma^{-1}$.
\item[Step 2]
Decide whether $\beta$ has a standard reduction system or not.
\item[Step 3]
If $\beta$ has no standard reduction system,
then return ``we cannot decide whether $\alpha$ is reducible'',
and halt.
\item[Step 4]
Find a standard reduction system, say $\C$, of $\beta$.
\item[Step 5]
Return ``$\gamma*\C$ is a reduction system of $\alpha$''.
\end{description}

Note that, from definitions,
$$[\alpha]^U\cap[\alpha]^U_\d =
\{ \beta\in [\alpha ]^S : \c^{\ell}(\beta)=\beta=\d^{m}(\beta)
\quad\mbox{for some $\ell, m\ge 1$} \}.
$$
This set is called the \emph{reduced super summit set},
and known to be nonempty~\cite{Lee00}.

Theorem~\ref{thm:main}~(i), (ii) and (iii) say that the above algorithm
finds a reduction system of a non-periodic reducible braid $\alpha$
if either $\infs(\alpha_\ext)>\infs(\alpha)$,
$\sups(\alpha_\ext)<\sups(\alpha)$,
or $\alpha$ is a split braid.
This implies that, roughly speaking,
if the outermost component $\alpha_\ext$ is simpler than the whole braid $\alpha$
up to conjugacy,
then we can find a reduction system of $\alpha$ from any element of
$[\alpha]^U\cap[\alpha]^U_\d$.

In Theorem~\ref{thm:main}, the conditions in (v) and (vi)
are weaker than those in (i) and (ii).
Because $\infs(\cdot)$ and $\sups(\cdot)$ are integer-valued,
Theorem~\ref{thm:INF}~(iii) implies the following.
\begin{itemize}
\item
If $\infs(\alpha_\ext)>\infs(\alpha)$,
then $\infs(\alpha_\ext)\ge\infs(\alpha)+1$ and, hence,
$$\INF(\alpha_\ext)\ge\infs(\alpha_\ext)\ge\infs(\alpha)+1>\INF(\alpha).$$
\item
If $\sups(\alpha_\ext)<\sups(\alpha)$,
then $\sups(\alpha_\ext)\le\sups(\alpha)-1$ and, hence,
$$\SUP(\alpha_\ext)\le\sups(\alpha_\ext)\le\sups(\alpha)-1<\SUP(\alpha).$$
\end{itemize}

Note that, for any $m\ne 0$, a braid $\alpha$ is reducible
if and only if $\alpha^m$ is reducible.
Therefore, in order to decide the reducibility of $\alpha$,
it suffices to decide the reducibility of $\alpha^m$
for an arbitrary $m\ne 0$.
If $\INF(\alpha_\ext)>\INF(\alpha)$ or $\SUP(\alpha_\ext)<\SUP(\alpha)$,
then the above algorithm,
applied to $\alpha^m$ for $1\le m< n(n-1)/2$,
finds a reduction system of $\alpha^m$
and, hence, decides the reducibility of $\alpha$.
Consequently, the non-periodic reducible braids
whose reducibility are not decidable by Theorem~\ref{thm:main}
are those with $\INF(\alpha_\ext)=\INF(\alpha)$
and $\SUP(\alpha_\ext)=\SUP(\alpha)$.

\medskip

We close this section with some examples.
From the examples, we can see that, in each statement of Theorem~\ref{thm:main},
the assertion does not hold if one of the conditions is weakened.

Example~\ref{ex:1} shows that Theorem~\ref{thm:main}~(i), (ii) and (iii) do not
hold for super summit sets.
Namely, there is a split braid
who satisfies the conditions (i) and (ii) but
whose super summit set contains an element without standard reduction system.

\begin{figure}
\begin{tabular}{ccc}
\includegraphics{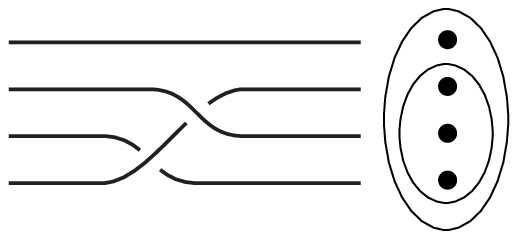} && \includegraphics{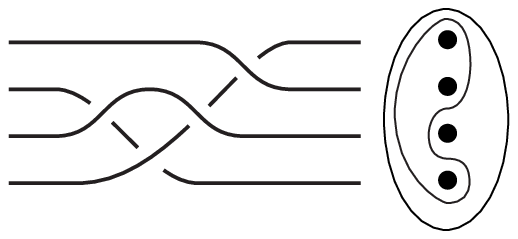}\\
\small (a) $\alpha=\sigma_1^{-1}\sigma_2 \in B_4$ &\qquad\quad&
\small (b) $\beta=\sigma_2^{-1}\sigma_1^{-1}\sigma_2\sigma_3$
\end{tabular}
\vskip -2mm \caption{$\alpha$ is a split braid. $\beta\in[\alpha]^S$
has no standard reduction system.}\label{fig:NonConvexSSS}
\end{figure}

\begin{example}\label{ex:1}
Let $\alpha=\sigma_1^{-1}\sigma_2 \in B_4$ and
$\beta=(\sigma_3\sigma_2)^{-1}\alpha (\sigma_3\sigma_2)
=\sigma_2^{-1}\sigma_1^{-1}\sigma_2\sigma_3$. (See
Figure~\ref{fig:NonConvexSSS}.)
Then $\alpha$ is a split braid with
$$0=\infs(\alpha_\ext) > \infs(\alpha)= -1 \quad\mbox{and}\quad
0= \sups(\alpha_\ext) < \sups(\alpha)= 1$$
and $\beta\in[\alpha]^S$, but $\beta$ has no standard reduction system.
\end{example}

Example~\ref{ex:2} shows the following.
\begin{itemize}
\item
Theorem~\ref{thm:main}~(i) and (ii) do not hold for
$\infs(\alpha_\ext)=\infs(\alpha)$ and $\sups(\alpha_\ext)=\sups(\alpha)$, respectively.
Namely, there is a non-periodic reducible braid
$\alpha$ with $\infs(\alpha_\ext)=\infs(\alpha)$ and
$\sups(\alpha_\ext)=\sups(\alpha)$ such that
the set $[\alpha]^U \cap [\alpha]^U_\d$ contains an element
without standard reduction system.

\item
For a non-periodic reducible braid $\alpha$ with periodic $\alpha_\ext$,
it is necessary to consider the ultra summit set $[\alpha^q]^U$
of some power of $\alpha$ in Theorem~\ref{thm:main}~(iv).
Namely, there is a non-periodic reducible $\alpha$ with
periodic $\alpha_\ext$ such that $[\alpha]^U$ contains
an element without standard reduction system.
\end{itemize}

\begin{figure}
\begin{tabular}{ccc}
\includegraphics{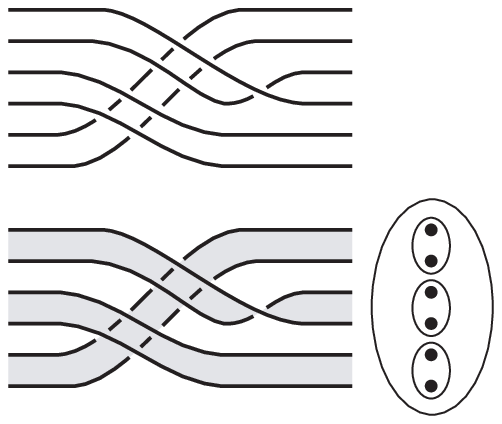} && \includegraphics{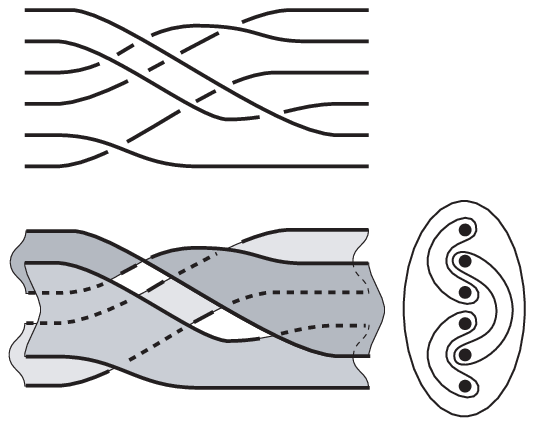} \\
(a) $\alpha=\sigma_2\sigma_1\sigma_3\sigma_2
\sigma_4\sigma_5\sigma_3\sigma_4\sigma_3$ &\qquad\qquad& (b) $\beta=
\sigma_4\sigma_1\sigma_3\sigma_2\sigma_4\sigma_5\sigma_4\sigma_3\sigma_2$
\end{tabular}
\vskip -2mm
\caption{The 6-braid $\alpha$ is non-periodic reducible
with $\infs(\alpha)=\infs(\alpha_\ext)$.
The braid $\beta$ belongs to $[\alpha]^U$, but $\beta$ has no
standard reduction system. }
\label{fig:NonConvexRSSS}
\end{figure}

\begin{example}\label{ex:2}
Consider the following 6-braids in Figure~\ref{fig:NonConvexRSSS}.
\begin{eqnarray*}
\alpha&=&\sigma_2\sigma_1\sigma_3\sigma_2
\sigma_4\sigma_5\sigma_3\sigma_4\sigma_3\\
\beta
&=& (\sigma_2\sigma_4^{-1})^{-1}\alpha(\sigma_2\sigma_4^{-1})
= \sigma_4\sigma_1\sigma_3\sigma_2\sigma_4\sigma_5\sigma_4\sigma_3\sigma_2
\end{eqnarray*}
Observe that $\alpha$ is a non-periodic reducible braid such that
$\alpha_\ext=\sigma_1\sigma_2$ is a periodic 3-braid.
Since $\alpha_\ext$, $\alpha$ and $\beta$ are all permutation braids,
we have
$$
\infs(\alpha)=0=\infs(\alpha_\ext);\quad
\sups(\alpha)=1=\sups(\alpha_\ext);\quad
\beta\in[\alpha]^U \cap [\alpha]^U_\d.
$$
It is easy to see that $\beta$ has no standard reduction system.
\end{example}

Example~\ref{ex:3} is due to Juan Gonz\'alez-Meneses and Bert Wiest.
The authors are very grateful to them for providing it.
It shows that
Theorem~\ref{thm:main}~(v) and (vi) do not hold for
$\INF(\alpha_\ext)=\INF(\alpha)$ and $\SUP(\alpha_\ext)=\SUP(\alpha)$,
respectively.
More precisely, there exist a non-periodic reducible braid $\alpha$
with $\INF(\alpha_\ext)=\INF(\alpha)$ and
$\SUP(\alpha_\ext)=\SUP(\alpha)$, and
an element $\beta$ such that, for each $q\ge 1$,
the power $\beta^q$ belongs to the set $[\alpha^q]^U\cap [\alpha^q]_\d^U$
but has no standard reduction system.

\begin{figure}
\begin{tabular}{ccc}
\includegraphics{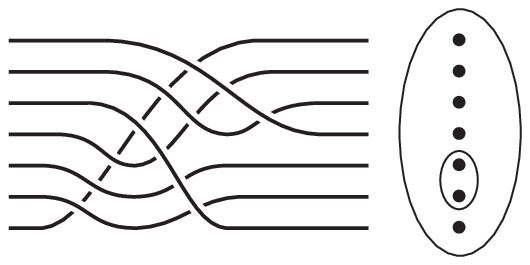} && \includegraphics{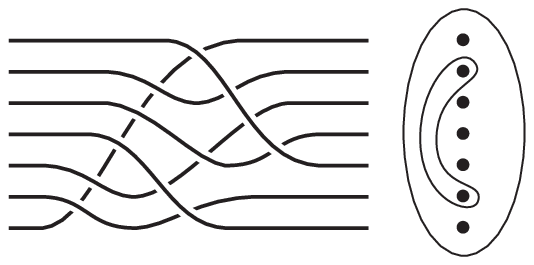} \\
(a) $\alpha=
\sigma_1\sigma_2\sigma_3\sigma_4\sigma_3\sigma_2\sigma_1
\sigma_5\sigma_4\sigma_6\sigma_5\sigma_4$ &\qquad\qquad& (b) $\beta=
\sigma_1\sigma_2\sigma_3\sigma_2\sigma_1\sigma_4\sigma_3
\sigma_5\sigma_6\sigma_5\sigma_4\sigma_3$
\end{tabular}
\vskip -2mm

\caption{The 7-braid $\alpha$ is non-periodic reducible
with $\INF(\alpha)=\INF(\alpha_\ext)=0$ and
$\SUP(\alpha_\ext)=\SUP(\alpha)=1$.
For all $q\ge 1$, the power $\beta^q$ belongs to
the set $[\alpha^q]^U\cap [\alpha^q]_\d^U$, but $\beta^q$ has no standard
reduction system.}
\label{fig:RedUSS}
\end{figure}

\begin{example}\label{ex:3}
Consider the following 7-braids in Figure~\ref{fig:RedUSS}.
\begin{eqnarray*}
\alpha&=&
\sigma_1\sigma_2\sigma_3\sigma_4\sigma_3\sigma_2\sigma_1
\sigma_5\sigma_4\sigma_6\sigma_5\sigma_4\\
\beta &=&
(\sigma_3\sigma_4\sigma_5)^{-1}
\alpha(\sigma_3\sigma_4\sigma_5)
= \sigma_1\sigma_2\sigma_3\sigma_2\sigma_1\sigma_4\sigma_3
\sigma_5\sigma_6\sigma_5\sigma_4\sigma_3
\end{eqnarray*}
Observe that
\begin{itemize}
\item[(i)]
both $\alpha$ and $\beta$ are permutation braids;
\item[(ii)]
$\alpha$ and $\beta$ are non-periodic reducible braids
with reduction systems as in Figure~\ref{fig:RedUSS};
\item[(iii)]
because $\alpha_\ext$ is pseudo-Anosov,
the curves in Figure~\ref{fig:RedUSS}~(a) and (b) are
the only reduction systems of $\alpha^q$ and $\beta^q$,
respectively, for all $q\ne 0$.
\end{itemize}
Let $B=\beta$.
(Throughout the paper, we have used capital letters
$A, B,\ldots$ to denote permutation braids.)
The starting set and finishing set of $B$ are
$$
S(B)=\{1,3,6\}\quad\mbox{and}\quad
F(B)=\{1,3,4,6\}.
$$
Since $S(B)\subset F(B)$,
the left normal form of $\beta^q$ is
$\Delta^0 \underbrace{BB\cdots B}_q$ for all $q\ge 1$.
In particular, for all $q\ge 1$,
$$
\c(\beta^q)=\beta^q,\quad
\d(\beta^q)=\beta^q,\quad
\inf(\beta^q)=0\quad\mbox{and}\quad
\sup(\beta^q)=q.
$$
Therefore, for all $q\ge 1$,
the power $\beta^q$ belongs to the set $[\alpha^q]^U\cap [\alpha^q]^U_\d$
and
\begin{eqnarray*}
\INF(\alpha)&=&\INF(\beta)=\lim_{q\to\infty}\inf(\beta^q)/q=0; \\
\SUP(\alpha)&=&\SUP(\beta)=\lim_{q\to\infty}\sup(\beta^q)/q=1.
\end{eqnarray*}
The outermost component $\alpha_\ext$ is obtained from $\alpha$
by deleting the second strand. Similarly to the above,
we can see that
$\INF(\alpha_\ext)=0=\INF(\alpha)$ and $\SUP(\alpha_\ext)=1=\SUP(\alpha)$.
\end{example}

\section{Proof of Proposition~\ref{thm:main-prop}}\label{sec:proof}
In this section, we prove Proposition~\ref{thm:main-prop}
that if $\alpha$ is a non-periodic reducible $n$-braid
with $\infs(\alpha_\ext)>\infs(\alpha)$,
then for any element $\beta$ of $[\alpha]^U$,
the cycling commutator $T_\beta$ is a split braid.

Throughout this section, the notation
$\St^\ext(\gamma)$ is used as an abbreviation for $\St(\Rext(\gamma))$,
the standardizer of the outermost component
of the canonical reduction system of the braid $\gamma$.
Therefore $\St^\ext(\gamma)$ consists of all positive braids $P$
such that $P*\Rext(\gamma)=\Rext(P\gamma P^{-1})$ is standard.
Recall that if $\gamma\in[\gamma]^U$
and $P$ is the $\le_R$-minimal element of $\St^\ext(\gamma)$,
then $P\gamma P^{-1}\in[\gamma]^U$ by Theorem~\ref{thm:StMain}.

\medskip
Let $\beta$ be an element of the ultra summit set $[\alpha]^U$.
Then $\c_0^m(\beta)=\beta$ for some $m\ge 1$.
For each $i=0,\ldots,m$, we define $n$-braids $A_i$,
$P_i$ and $\gamma^{(i)}$ as follows (see Figure~\ref{fig:comm}):
\begin{itemize}
\item $A_i$ is the $\le_R$-minimal element of
$\{P\in B_n^+: \inf(P\c_0^i(\beta))>\inf(\c_0^i(\beta))\}$;
\item $P_i$ is the $\le_R$-minimal element of $\St^\ext(\c_0^i(\beta))$;
\item $\gamma^{(i)}=P_i \c_0^i(\beta) P_i^{-1}$.
\end{itemize}
\begin{figure}
$$
\begin{CD}
\beta @>A_0>> \c_0(\beta) @>A_1>> \c_0^2(\beta) @>A_2>>
\quad\cdots\quad  @>A_{m-1}>> \c_0^m(\beta)=\beta\\
@VV{P_0}V    @VV{P_1}V    @VV{P_2}V    @.      @VV{P_m=P_0}V\\
\gamma^{(0)} @>B_0>> \gamma^{(1)} @>B_1>> \gamma^{(2)} @>B_2>>
\quad\cdots\quad  @>B_{m-1}>> \gamma^{(m)}=\gamma^{(0)}\\
\end{CD}
$$
\vskip -3mm
\caption{``$\alpha\stackrel{A}{\longrightarrow}\beta$''
denotes $\beta=A\alpha A^{-1}$.}
\label{fig:comm}
\end{figure}
Then, for each $i=0,\ldots, m-1$,
\begin{itemize}
\item
$A_i$ is a permutation braid with
$\c_0^{i+1}(\beta)=A_i \c_0^i(\beta) A_i^{-1}$ by Lemma~\ref{thm:duality}~(i);

\item
$\Rext(\gamma^{(i)})$ is standard because
$\Rext(\gamma^{(i)})=\Rext(P_i\c_0^i(\beta) P_i^{-1})=P_i*\Rext(\c_0^i(\beta))$
and $P_i\in\St^\ext(\c_0^i(\beta))$;

\item
$\gamma^{(i)}$ belongs to $[\alpha]^U$ by Theorem~\ref{thm:StMain}.
\end{itemize}

\begin{lemma}\label{lemma:1}
For $i=0,\ldots,m-1$, there exists a permutation braid $B_i$
such that $B_iP_i=P_{i+1}A_i$ and $\gamma^{(i+1)}=B_i\gamma^{(i)}B_i^{-1}$.
\end{lemma}

\begin{figure}
$
\xymatrix{ \c_0^i(\beta) \ar[r]^{A_i} \ar[d]^{P_i}
&  \c_0^{i+1}(\beta)  \ar[d]^{P_{i+1}} \ar@/^1em/[ddr]^{P_{i+1}'} \\
\gamma^{(i)} \ar[r]^{B_i} \ar@/_1em/[drr]^{B_i'}
&  \gamma^{(i+1)}  \ar[dr]^{B_i''} \\
& & \c_0(\gamma^{(i)}) }
$
\vskip -2mm
\caption{``$\alpha\stackrel{A}{\longrightarrow}\beta$''
denotes $\beta=A\alpha A^{-1}$.}
\label{fig:comm1}
\end{figure}

\begin{proof}
(See Figure~\ref{fig:comm1}.)
Let $B_i'$ be the $\le_R$-minimal element of
$\{P\in B_n^+:\inf(P\gamma^{(i)})>\inf(\gamma^{(i)})\}$.
Then $B_i'$ is a permutation braid by Lemma~\ref{thm:duality},
and
$$
\inf(B_i'\gamma^{(i)})>\inf(\gamma^{(i)})
\quad\mbox{and}\quad
\c_0(\gamma^{(i)})=B_i'\gamma^{(i)}B_i'^{-1}.
$$
Since both $\gamma^{(i)}$ and $\c_0^i(\beta)$ belong to $[\alpha]^U$,
we have $\inf(\gamma^{(i)})=\inf(\c_0^i(\beta))=\infs(\alpha)$.
Since
$$
\inf(B_i'P_i \c_0^i(\beta))=\inf(B_i'\gamma^{(i)}P_i)
\ge\inf(B_i'\gamma^{(i)})>\inf(\gamma^{(i)})=\inf(\c_0^i(\beta)),
$$
$B_i'P_i$ belongs to the set
$\{P\in B_n^+: \inf(P\c_0^i(\beta))>\inf(\c_0^i(\beta))\}$.
Since $A_i$ is the $\le_R$-minimal element of this set,
we have $A_i\le_R B_i'P_i$, and hence
\begin{equation}\label{eqn:8-1}
B_i'P_i=P_{i+1}'A_i
\end{equation}
for some $P_{i+1}'\in B_n^+$.
Note that
$$
P_{i+1}'\c_0^{i+1}(\beta)P_{i+1}'^{-1}
= P_{i+1}'A_i \c_0^i(\beta) A_i^{-1} P_{i+1}'^{-1}
= B_i'P_i \c_0^i(\beta) P_i^{-1} B_i'^{-1}
= B_i'\gamma^{(i)}B_i'^{-1}
=\c_0(\gamma^{(i)}).
$$
Since
$\Rext(\c_0(\gamma^{(i)}))$ is standard by Lemma~\ref{thm:standardSSS},
$P_{i+1}'$ belongs to $\St^\ext(\c_0^{i+1}(\beta))$.
Since $P_{i+1}$ is the $\le_R$-minimal element of $\St^\ext(\c_0^{i+1}(\beta))$,
we have $P_{i+1}\le_R P_{i+1}'$.
Therefore,
\begin{equation}\label{eqn:8-2}
P_{i+1}'=B_i''P_{i+1}
\end{equation}
for some $B_i''\in B_{n}^+$.
Observe that
$$
P_{i+1}A_i \c_0^{i}(\beta) A_i^{-1}P_{i+1}^{-1}
=P_{i+1}  \c_0^{i+1}(\beta) P_{i+1}^{-1}
=\gamma^{(i+1)}.
$$
Since $\R_\ext(\gamma^{(i+1)})$ is standard,
$P_{i+1}A_i$ belongs to $\St^\ext(\c_0^{i}(\beta))$.
Since $P_i$ is the $\le_R$-minimal element of
$\St^\ext(\c_0^{i}(\beta))$, we have $P_i\le_R P_{i+1}A_i$.
Therefore
\begin{equation}\label{eqn:8-3}
P_{i+1}A_i=B_iP_i
\end{equation}
for some $B_i\in B_n^+$.
It is obvious that $\gamma^{(i+1)}=B_i\gamma^{(i)}B_i^{-1}$.
From (\ref{eqn:8-1}), (\ref{eqn:8-2}) and (\ref{eqn:8-3}),
$$B_i'P_i=P_{i+1}'A_i=B_i''P_{i+1}A_i=B_i''B_iP_i.$$
Therefore $B_i'=B_i''B_i$.
Since $B_i'$ is a permutation braid and $B_i\le_R B_i'$,
the positive braid $B_i$ is a permutation braid as desired.
\end{proof}

Let $\Rext(\gamma^{(0)})=\C_\n$ for a composition
$\n=(n_1,\ldots,n_k)$ of $n$.
Let $\Delta_i$ be the fundamental braid of $B_{n_i}$.

\begin{lemma}\label{lemma:2}
For $i=0,\ldots,m-1$,
$\Rext(\gamma^{(i)})=\C_\n$ and
$B_i\le_R(\Delta_{1}\oplus\cdots\oplus\Delta_{k})$.
\end{lemma}

\begin{proof}
Using induction on $i$, it suffices to show the following:

\medskip
\begin{quote}\em
If\/ $\Rext(\gamma^{(i)})=\C_\n$, then
$B_i\le_R(\Delta_{1}\oplus\cdots\oplus\Delta_{k})$
and $\Rext(\gamma^{(i+1)})=\C_\n$.
\end{quote}
\medskip

Suppose $\Rext(\gamma^{(i)})=\C_\n$.
By Lemma~\ref{thm:decom}~(ii) and (iv),
$$
\gamma^{(i)}
=(\gamma_1\oplus\cdots\oplus\gamma_k)\myangle{\gamma_0}_\n,
$$
where $\gamma_0=\gamma^{(i)}{}_\ext\in B_k$
and $\gamma_j\in B_{n_j}$ for $j=1,\ldots,k$.
Since $\infs(\alpha_\ext)>\infs(\alpha)$ (from the hypothesis)
and $\gamma^{(i)}\in[\alpha]^U$,
we have $\inf(\gamma^{(i)}{}_\ext)>\inf(\gamma^{(i)})$
by Lemma~\ref{thm:inf_out}. By Lemma ~\ref{thm:infofredbr},
$$\inf(\gamma_0)=\inf(\gamma^{(i)}{}_\ext)>\inf(\gamma^{(i)})
=\min\{\inf(\gamma_i):i=0,\ldots,k,\ \brindex(\gamma_i)\ge 2 \}.
$$
Therefore $\inf(\gamma_0)>\inf(\gamma_j)$
for some $j\ge 1$ with $\brindex(\gamma_j)\ge 2$, and
\begin{eqnarray*}
\lefteqn{\inf((\Delta_{1}\oplus\cdots\oplus\Delta_{k})\gamma^{(i)})
=\inf((\Delta_1\gamma_1\oplus\cdots\oplus\Delta_k\gamma_k)\myangle{\gamma_0}_\n)}\\
&=& \min\ (\{\inf(\Delta_j\gamma_j):j=1,\ldots,k, \brindex(\gamma_j)\ge 2\}
\cup\{\inf(\gamma_0)\})\\
&=& \min\ (\{\inf(\gamma_j)+1:j=1,\ldots,k, \brindex(\gamma_j)\ge 2\}
 \cup\{\inf(\gamma_0)\})\\
&>&\inf(\gamma^{(i)}).
\end{eqnarray*}
So $(\Delta_{1}\oplus\cdots\oplus\Delta_{k})\in
\{P\in B_n^+:\inf(P\gamma^{(i)})>\inf(\gamma^{(i)})\}$.
Recall, from the proof of Lemma~\ref{lemma:1}, that
$B_i\le_R B_i'$, where $B_i'$ is the $\le_R$-minimal element of
$\{P\in B_n^+:\inf(P\gamma^{(i)})>\inf(\gamma^{(i)})\}$.
Therefore,
$$B_i\le_R B_i'\le_R(\Delta_1\oplus\cdots\oplus\Delta_k)$$
as desired.
This implies that $B_i$ has the decomposition
$B_i=(B_{i,1}\oplus\cdots\oplus B_{i,k})$
for some permutation $n_j$-braid $B_{i,j}$'s.
By Lemma~\ref{thm:decom}~(ii), $B_i*\C_\n=\C_\n$.
Therefore,
$\Rext(\gamma^{(i+1)})=\Rext(B_i\gamma^{(i)}B_i^{-1})
=B_i*\Rext(\gamma^{(i)})=B_i*\C_\n=\C_\n$.
\end{proof}

Let $S=B_{m-1}\cdots B_0$.
Then $S$ is a split braid by Lemma~\ref{lemma:2}.
Note that the cycling commutator of $\beta$ is
$T_{\beta}=A_{m-1}\cdots A_0$.
Since $P_0^{-1}SP_0=T_{\beta}$ by Lemma~\ref{lemma:1},
$T_{\beta}$ is a split braid and the proof is completed.

\end{document}